\documentclass{article}
\usepackage{fullpage}
\usepackage{amsmath,amsthm,mathrsfs, verbatim,color,lscape,xspace}
\allowdisplaybreaks
\usepackage{listings}
\usepackage{graphicx}
\usepackage{secdot}
\usepackage{amsmath, dsfont}
\usepackage{amsfonts}
\usepackage{amssymb}
\usepackage{enumerate}
\usepackage{xspace}
\usepackage{enumitem}
\usepackage{ifpdf}
\usepackage{lscape}
\usepackage[tight]{subfigure}
\usepackage{hyperref}      %Use only when compiling pdf's, otherwise it creates errors
\usepackage{xcolor}\hypersetup{pdftex,breaklinks,citebordercolor=green,linkbordercolor=magenta}

\newcommand{\cordisloss}[1][t]{\ensuremath{Loss_{d,\epsilon}(x,#1)}}         %for the Loss corrected discard approximation
\newcommand{\loss}[1][t]{\ensuremath{Loss_{\epsilon}(x,#1)}}         %for the loss function
\newcommand{\ltptc}[2][s]{\ensuremath{\beta^{#2}(#1)}\xspace}        %for the L-T of the phase-type claim sizes
\newcommand{\lttime}[1][s]{\ensuremath{\tau\left(#1\right)}\xspace}        %for the L-T of variable T

\newcommand{\erfct}[1]{\frac{2}{\sqrt\pi}\int_{#1}^\infty e^{-x^2}dx}

\newcommand{\lt}[1]{\ensuremath{#1}}          %for L.T. of f, I have to type \LT{f}, argument not optional

\newcommand{\com}[1]{\ensuremath{\overline{#1}}}

\newcommand{\disf}[2][u]{\ensuremath{{#2}(#1)}}     %Distribution  function of f with default parameter t,
\newcommand{\df}[2][u]{\ensuremath{{#2}(#1)}}       %it gives the density function of f with default parameter t,
\newcommand{\con}[3][n]{\ensuremath{{#2}^{*#1}_{#3}}}            %[#1]=n'th convolution,{#2}=distribution,{#3}=index
\newcommand{\e}[1][]{\ensuremath{\mathbb{E}{#1}}\xspace}    %for the expectation
\newcommand{\pr}[1][]{\ensuremath{\mathbb{P}{#1}}\xspace}   %for the probability
\newcommand{\ruinmix}[1][u]{\ensuremath{\psi_\epsilon(#1)}}                %for the ruin probability of the mixed model (notation with indices)
\newcommand{\ruindis}[1][u]{\ensuremath{\psi^\bullet_\epsilon(#1)}}        %for the ruin probability of the discard model
\newcommand{\ruinrep}[1][u]{\ensuremath{\psi_0(#1)}}                      %for the ruin probability of the replace model

\newcommand{\cordis}[1][u]{\ensuremath{\tilde{\psi}_{d,\epsilon}(#1)}}         %for the corrected discard approximation
\newcommand{\correp}[1][u]{\ensuremath{\tilde{\psi}_{r,\epsilon}(#1)}}         %for the corrected replace approximation

\newcommand{\ltruinmix}[1][s]{\ensuremath{m_\epsilon(#1)}\xspace}              %for the L-T of the mixed ruin probability
\newcommand{\ltruindis}[1][s]{\ensuremath{m^\bullet_\epsilon(#1)}\xspace}              %for the L-T of the discard ruin probability
\newcommand{\ltruinrep}[1][s]{\ensuremath{m_0(#1)}\xspace}              %for the L-T of the replace ruin probability

\newcommand{\ltegc}[1][s]{\ensuremath{\upsilon_\epsilon^e(#1)}\xspace}        %for the L-T of the stationary-excess general claim sizes
\newcommand{\ltgc}[1][s]{\ensuremath{\upsilon_\epsilon(#1)}\xspace}        %for the L-T of the general claim sizes
\newcommand{\lteptc}[1][s]{\ensuremath{\beta^e(#1)}\xspace}        %for the L-T of the stationary-excess phase-type claim sizes
\newcommand{\ltehtc}[1][s]{\ensuremath{\gamma^e(#1)}\xspace}        %for the L-T of the stationary-excess heavy-tailed claim
\newcommand{\lthtc}[1][s]{\ensuremath{\gamma(#1)}\xspace}        %for the L-T of the stationary-excess heavy-tailed claim

\newcommand{\ptc}[1][]{\ensuremath{B_{#1}}\xspace}                    %for the phase-type service times
\newcommand{\eptc}[1][]{\ensuremath{B^e_{#1}}\xspace}                 %for the stationary-excess phase-type service times
\newcommand{\htc}[1][]{\ensuremath{C_{#1}}\xspace}                    %for the heavy-tailed service times
\newcommand{\ehtc}[1][]{\ensuremath{C^e_{#1}}\xspace}                 %for the stationary-excess heavy-tailed service times
\newcommand{\mixc}[1][]{\ensuremath{U_{\epsilon #1}}\xspace}          %for general mixed service times
\newcommand{\emixc}[1][]{\ensuremath{U^e_{\epsilon #1}}\xspace}       %for stationary-excess general mixed service times
\newcommand{\mixsup}[1][]{\ensuremath{M_{\epsilon #1}}\xspace}        %for the supremum of the mixed claim surplus process
\newcommand{\dissup}[1][]{\ensuremath{M^\bullet_{\epsilon #1}}\xspace}        %for the supremum of the discard claim surplus process
\newcommand{\repsup}[1][]{\ensuremath{M_{0 #1}}\xspace}        %for the supremum of the replace claim surplus process

\newtheorem{lemma}{Lemma}

\newtheorem{definition}{Definition}
\newtheorem{theorem}{Theorem}
\theoremstyle{definition}
\newtheorem{remark}{Remark}
\newtheorem{property}{Property}

\newcommand{\footnoteremember}[2]
{
   \newcounter{#1}\footnote{#2}\setcounter{#1}{\value{footnote}}
}
\newcommand{\footnoterecall}[1]
{
   \footnotemark[\value{#1}]
}

\begin{document}
\title{Corrected phase-type approximations of heavy-tailed risk models using perturbation analysis}
\author{
     E. Vatamidou\footnoteremember{TU/eEURANDOM}{\textsc{Eurandom} and Department of Mathematics \& Computer Science, Eindhoven University of Technology, P.O. Box 513, 5600 MB Eindhoven, The Netherlands}\\
     \small \texttt{e.vatamidou@tue.nl}\\
     \and
     I.J.B.F. Adan\footnoterecall{TU/eEURANDOM}\footnoteremember{MechEng}{Department of Mechanical Engineering, Eindhoven University of Technology, P.O. Box 513, 5600 MB Eindhoven, The Netherlands}\\
     \small \texttt{i.j.b.f.adan@tue.nl}\\
     \and
     M. Vlasiou\footnoterecall{TU/eEURANDOM}\footnoteremember{CWI}{Centrum Wiskunde \& Informatica (CWI), P.O. Box 94079, 1090 GB Amsterdam, The Netherlands}\\
     \small \texttt{m.vlasiou@tue.nl}
     \and
     B. Zwart\footnoterecall{TU/eEURANDOM}\footnoterecall{CWI}\\
     \small \texttt{Bert.Zwart@cwi.nl}\\
}
\maketitle

\begin{abstract}
  Numerical evaluation of performance measures in heavy-tailed risk models is an important and challenging problem. In this paper, we construct very accurate approximations of such performance measures that provide small absolute and relative errors. Motivated by statistical analysis, we assume that the claim sizes are a mixture of a phase-type and a heavy-tailed distribution and with the aid of perturbation analysis we derive a series expansion for the performance measure under consideration. Our proposed approximations consist of the first two terms of this series expansion, where the first term is a phase-type approximation of our measure. We refer to our approximations collectively as {\it corrected phase-type approximations}. We show that the corrected phase-type approximations exhibit a nice behavior both in finite and infinite time horizon, and we check their accuracy through numerical experiments.
\end{abstract}

\section{Introduction}
The evaluation of performance measures of risk models is an important problem that has been widely studied in the literature \cite{asmussen-RP,kluppelberg04,kyprianou-ILFLP}. Under the presence of heavy-tailed claim sizes, these evaluations become more challenging and sometimes even problematic \cite{Ahn12,asmussen05}. In such cases, it is necessary to construct approximations for the quantity under consideration. In this paper, we develop a new method to construct reliable approximations for performance measures of heavy-tailed risk models. We use the classical risk model (perhaps outdated, but very well studied) as a context and vehicle to demonstrate our key ideas, which we expect to have a much wider applicability in insurance. We show that our approximations have a provably small absolute error, independent of the initial capital, and a small relative error. As additional test of performance we also consider the finite horizon aggregate loss model.

There are three main directions for approximating ruin probabilities: phase-type approximations, asymptotic approximations and error bounds. When the claim sizes follow some light-tailed distribution, a natural approach to provide approximations for the ruin probability with high accuracy is by approximating the claim size distribution with a phase-type one \cite{feldmann98,sasaki04,starobinski00}. We refer to these methods as {\it phase-type approximations}, because the approximate ruin probability has a phase-type representation \cite{asmussen92a,ramaswami90}. However, in many financial applications, an appropriate way to model claim sizes is by using heavy-tailed distributions \cite{asmussen-APQ,embrechts-MEE,rolski-SPIF}. In these cases, the exponential decay of phase-type approximations gives a big relative error at the tail  and the evaluation of the ruin probability becomes more complicated.

When the claim size distribution belongs to the class of subexponential distributions \cite{teugels75}, which is a special case of heavy-tailed distributions, asymptotic approximations are available \cite{bahr75,borovkov92,embrechts82,olvera11,pakes75}. The main disadvantage of such approximations is that they provide a good fit only at the tail of the ruin probability, especially for small safety loading. Another stream of research focuses on corrected diffusion approximations for the ruin probability \cite{blanchet10,silvestrov-LTRSSP}. A disadvantage of such asymptotic techniques is the requirement of finite higher moments for the claim size distribution.

Finally, results on error bounds \cite{kalashnikov02,vatamidou12} indicate that such bounds are rather pessimistic, especially in terms of relative errors, and in case of small safety loading. There exist also bounds with the correct tail behavior under subexponential claims \cite{kalashnikov02a,korshunov11}, but these bounds are only accurate at the tail. A conclusion that can be safely drawn from all the above is that, although the literature is abundant with approximations for the ruin probability in the case of light-tailed claim sizes, accurate approximations for the ruin probability in the case of heavy-tailed claim sizes are still an open topic.

Besides the ruin probability, a very popular tool in real-world applications to measure the operational risk is the Value at Risk (VaR) \cite{embrechts-QRM}. For a given portfolio, a VaR with a probability level $\alpha$ and fixed time horizon is defined as the threshold value such that the loss on the portfolio over the given time horizon exceeds this value with probability $1-\alpha$. It is of interest to quantify the operational risk through the statistical analysis of operational loss data \cite{embrechts02,klugman-LM} and provide error bounds for the aggregate loss probability \cite{cox08}. Similarly to the ruin probability, things become more complicated under the presence of heavy-tailed data \cite{embrechts-MEE}.

In this paper, we develop approximations for ruin probabilities and total losses under heavy-tailed claims that combine desirable characteristics of all three main approximation directions. First, our approximations maintain the computational tractability of phase-type approximations. Additionally, they capture the correct tail behavior, which so far could only be captured by asymptotic approximations, and they have the advantage that finite higher-order moments are not required for the claim sizes. Last, they provide a provably small absolute error, independent of the initial capital, and a small relative error.

The idea of our approach stems from fitting procedures of the claim size distribution to data. Heavy-tailed statistical analysis suggests that for a sample with size $n$ only a small fraction ($k_n/n \rightarrow 0$) of the upper-order statistics is relevant for estimating tail probabilities \cite{davis84,hill75,resnick-HTP}. More information about the optimal choice of the $k_n$th upper order statistic can be found in \cite{haeusler85}. The remaining data set may be used to fit the bulk of the distribution. Since the class of phase-type distributions is dense in the class of all positive definite probability distributions \cite{asmussen-APQ}, a natural choice is to fit a phase-type distribution to the remaining data set \cite{asmussen96b}. As a result, a mixture model for the claim size distribution is a natural assumption. Thus, our key idea is to use a mixture model for the claim size distribution in order to construct approximations of the ruin probability that combine the best elements of phase-type and asymptotic approximations.

We now sketch how to derive our approximations when the claim size distribution is a mixture of a phase-type distribution and a heavy-tailed one. Interpreting the heavy-tailed term of the claim size distribution in the mixture model as perturbation of the phase-type one and using perturbation theory, we can find the ruin probability (total loss) as a complete series expansion. The first term of the expansion is the phase-type approximation of the ruin probability (total loss) that occurs when we ``remove" the heavy-tailed claim sizes from the system, either by discarding them or by replacing them with phase-type ones. We consider the model that appears when all heavy-tailed claims are removed as the ``base" model. Due to the two different approaches of removing the heavy-tailed claim sizes, the ruin probability (total loss) connects to two different base models and consequently to two different series expansions.

We show that adding the second term of the respective series expansions is sufficient to construct improved approximations, compared to their phase-type counterparts, the {\it discard} and the {\it replace approximations}, respectively. Since the second term of each series expansion works as a correction to its respective phase-type approximation, motivated by the terminology {\it corrected heavy traffic approximations} \cite{asmussen-APQ}, we refer to our approximations as {\it corrected phase-type approximations}. Therefore, in this paper, we propose the {\it corrected discard approximation} and the {\it corrected replace approximation}. Both approximations have appealing properties: the corrected replace approximation tends to give better numerical estimates, while the corrected discard approximation is simpler and yields guaranteed upper and lower bounds. Last, we provide the form of the {\it corrected phase-type approximations} for the aggregate loss over a fixed time period, and we show that they have the same appealing properties also for finite time.

Within risk theory, some attention has been given to perturbed risk models; see \cite{schmidli-PRP} for a review and the recent paper of \cite{huzak04}. However, the term ``perturbation" in this area is used to denote the superposition of two risk processes. Contrary to other asymptotic techniques that use perturbation analysis to approximate the ruin probability \cite{blanchet10,silvestrov-LTRSSP}, our approach is different; we apply perturbation to the claim sizes rather than the arrival rate.

The connection between ruin probabilities and the stationary waiting probability $\pr (W_q>u)$ of a G/G/1 queue, where service times in the queueing model correspond to the random claim sizes, is well known \cite{asmussen-APQ,asmussen-RP}. Thus, the corrected approximations can also be used to estimate the waiting time distribution of the above mentioned queue. Finally, since the reserve process of the classical risk model is a basic building block of any L\'{e}vy process \cite{kluppelberg04,kyprianou-ILFLP}, and due to the connection of ruin probabilities with scale functions \cite{Ahn12,biffis10}, we expect that our technique is widely applicable to more general risk processes.

The rest of the paper is organized as follows. In Section~\ref{S.Power series expansion}, we introduce the model and we derive two series expansions for the ruin probability. From these series expansions we deduce approximations for the ruin probability, in Section~\ref{S.Approximations}, and we study their basic properties. In Section~\ref{S.Example}, we find the exact formula of the ruin probability for a specific mixture model and we study the extent of the achieved improvement when we compare our approximations with phase-type approximations of their related base model. In Section~\ref{S.VaR}, we provide corrected phase-type approximations of the aggregate loss in finite time and we show through a numerical study that our approximations give excellent VaR estimates. Finally, in the Appendix, we give all the proofs.

\section{Series expansions of the ruin probability}\label{S.Power series expansion}
As proof of concept, we apply our technique to the classical Cram\'{e}r-Lundberg risk model \cite{asmussen-RP,prabhu61}. In this model, we assume that premiums flow in at a rate 1 per unit time and claims arrive according to a Poisson process $\{\df[t]{N_\epsilon}\}_{t \geq 0}$ with rate $\lambda$, where $\epsilon \in [0,1]$ is a parameter to be explained soon. The claim sizes $\mixc[,i]\stackrel{d}{=}\mixc$ are i.i.d.\ with common distribution $G_\epsilon$ and independent of $\{\df[t]{N_\epsilon}\}$. Motivated by statistical analysis, which proposes that only a small fraction of the upper-order statistics is relevant for estimating tail probabilities, we consider that an arbitrary claim size \mixc\ is phase-type \cite{neuts-MGSSM} with probability $1-\epsilon$ and heavy-tailed \cite{rolski-SPIF} with probability $\epsilon$, where $\epsilon \rightarrow 0$. In the forthcoming analysis, we use as general rule that all parameters depending on $\epsilon$ bear a subscript with the same letter. We assume that the phase-type claim sizes $\ptc[i] \stackrel{d}{=} \ptc$ and the heavy-tailed claim sizes $\htc[i] \stackrel{d}{=} \htc$ have both finite means, $\e{\ptc}$ and $\e{\htc}$, respectively. If $u$ is the initial capital, our risk reserve process has the form
\begin{equation*}
  \df[t]{R_\epsilon} = u+t-\sum_{i=1}^{\df[t]{N_\epsilon}}\mixc[,i].
\end{equation*}

Using this model, we first examine in Sections~\ref{S.Power series expansion}--\ref{S.Example} the ruin probability in infinite time horizon, and later on, in Section~\ref{S.VaR}, we move to finite time horizon and we examine the aggregate loss.

When calculating ruin probabilities, for mathematical purposes, it is more convenient to work with the claim surplus process $\df[t]{S_\epsilon} =u - \df[t]{R_\epsilon}$. The probability \ruinmix\ of ultimate ruin is the probability that the reserve ever drops below zero or equivalently the probability that the maximum $\mixsup = \sup_{0 \leq t < \infty} \df[t]{S_\epsilon}$ ever exceeds $u$; i.e.
\begin{equation*}
  \ruinmix = \pr(\mixsup > u).
\end{equation*}

For a distribution $F$ we use the notation $\con{F}{}$ for its $n$th convolution and $\com{F}$ for its complementary cumulative distribution $1-F$. Moreover, if $F$ has a finite mean $\mu_F$, then we define its stationary excess distribution as
\begin{equation*}
  \disf{F^e} = \frac{1}{\mu_F} \int_0^u \disf[t]{\com{F}}dt.
\end{equation*}
In addition, the r.v.\ with distribution $F^e$ bears also a superscript with the letter ``{\it e}".

When the average amount of claim per unit time $\rho_\epsilon = \lambda \e{\mixc}$ is strictly smaller than 1, the well-known Pollaczek-Khinchine formula \cite{asmussen-RP} can be used for the evaluation of the ruin probability. We, namely, have that
\begin{equation}\label{e.PK perturbed ruin probability}
  1 - \ruinmix = (1-\rho_\epsilon)\sum_{n=0}^\infty \rho_\epsilon^n \df{{\con{\left(G^e_\epsilon\right)}{}}},
\end{equation}
where $G^e_\epsilon$ is the distribution of the stationary excess claim sizes \emixc. The infinite sum of convolutions at the right-hand side of \eqref{e.PK perturbed ruin probability} makes the evaluation of \ruinmix\ difficult or even impossible for our mixture model. For this reason, one typically resorts to Laplace transforms. We use the notation \ltegc, \lteptc\ and \ltehtc\ for the Laplace transforms of the stationary excess claim sizes $\emixc[,i]\stackrel{d}{=}\emixc$, $\eptc[i]\stackrel{d}{=}\eptc$ and $\ehtc[i]\stackrel{d}{=}\ehtc$, respectively. Moreover, we set $\delta=\lambda \e{\ptc}$ and $\theta = \lambda \e{\htc}$, which means that the phase-type claims are responsible for average claim $(1-\epsilon)\delta$ per unit time and the heavy-tailed claims are responsible for average claim $\epsilon\theta$ per unit time. Using this notation, we obtain $\rho_\epsilon = (1-\epsilon)\delta + \epsilon\theta$. In terms of Laplace transforms, the Pollaczek-Khinchine formula can be written now as:
\begin{equation}\label{e.LT PK perturbed ruin probability}
    \ltruinmix := \e e^{-s \mixsup} \
                = (1-\rho_\epsilon)\sum_{n=0}^\infty\rho_\epsilon^n \left(\ltegc\right)^n
                = \frac{1-\rho_\epsilon}{1 - \rho_\epsilon \ltegc}
                = \frac{1-(1-\epsilon)\delta-\epsilon\theta}{1-(1-\epsilon)\delta \lteptc - \epsilon\theta \ltehtc}.
\end{equation}

Applying Laplace inversion to \eqref{e.LT PK perturbed ruin probability} to find \ruinmix\ is difficult \cite{abate99b} or even impossible, because the heavy-tailed component \ltehtc\ oftentimes does not have an analytic closed form. To overcome this difficulty, a phase-type approximation would suggest to ``remove" the heavy-tailed claim sizes and find an explicit phase-type representation for the ruin probability of the resulting simpler model, which we use as base model for our analysis. In broad terms, we view the heavy-tailed claim sizes as perturbation of the phase-type claim sizes and we interpret $\epsilon$ as the perturbation parameter. With the aid of perturbation analysis, we find the ruin probability of our mixture model as a complete series expansion with first term the phase-type approximation that results from its base model.

As mentioned in the introduction, we remove the heavy-tailed claims either by discarding them or by replacing them with phase-type ones. Therefore, the ruin probability \ruinmix\ connects to two different base models and, consequently, it has two different series expansions, the discard and the replace expansions. We first derive the discard series expansion. From a mathematical point of view, when we discard the heavy-tailed claim sizes, we simply consider that $\disf[x]{G_\epsilon} = (1-\epsilon) \pr\left(\ptc \leq x\right) + \epsilon$, $x\geq 0$. This base model, for which the claim size distribution has an atom at zero, is equivalent to the compound Poisson risk model in which claims arrive with rate $(1-\epsilon)\lambda$ and follow the distribution of \ptc. We denote by \dissup\ the supremum of its corresponding claim surplus process. Thus, the Pollaczek-Khinchine formula for this base model takes the form
\begin{equation}\label{e.LT PK discard ruin probability}
  \ltruindis  := \e e^{-s \dissup} = \frac{1-(1-\epsilon)\delta}{1-(1-\epsilon)\delta \lteptc}.
\end{equation}
We denote by \ruindis\ the discard phase-type approximation of \ruinmix\ that appears when we apply Laplace inversion to the above formula. For this base model, the series expansion of \ruinmix\ can be found in the following theorem.

\begin{theorem}\label{t.discard expansion}
  \textup{\textbf{Discard expansion.}} If \ruindis\ is the phase-type approximation of the exact ruin probability \ruinmix\ that occurs when we discard the heavy-tailed claim sizes and $\dissup[,i]\stackrel{d}{=}\dissup$, a series expansion of the exact ruin probability is given by
  \begin{align*}
    \ruinmix   = \ruindis
                        &+ \sum_{n=1}^\infty \left(\frac{\epsilon\theta}{1-\delta+\epsilon\delta}\right)^n \big(\disf{L_{\epsilon,n}} - \disf{L_{\epsilon,n-1}} \big),
  \end{align*}
  where $\disf{L_{\epsilon,n}} = \pr(\dissup[,0]+\dissup[,1]+\dots+\dissup[,n]+\ehtc[1]+\dots+\ehtc[n] > u)$ and $\disf{L_{\epsilon,0}} = \pr(\dissup[,0] > u) = \ruindis$. A necessary and sufficient condition for the convergence of the discard series expansion for all values of $u$ is $\lvert \epsilon\theta \rvert < \lvert 1-\delta+\epsilon\delta \rvert$.
\end{theorem}

To find the replace series expansion, observe that the action of replacing the heavy-tailed claim sizes with phase-type ones translates into $\epsilon=0$. For this base model, the Pollaczeck-Khinchine formula takes the form
\begin{equation}\label{e.LT PK replace ruin probability}
  \ltruinrep  := \e e^{-s \repsup} = \frac{1-\delta}{1-\delta \lteptc},
\end{equation}
where $\repsup = \mixsup[|\epsilon=0]$. Laplace inversion of \ltruinrep\ gives the phase-type approximation \ruinrep\ of the ruin probability \ruinmix. The series expansion of \ruinmix\ in this case is given below.

\begin{theorem}\label{t.replace expansion}
  \textup{\textbf{Replace expansion.}} If \ruinrep\ is the phase-type approximation of the exact ruin probability \ruinmix\ that occurs when we replace the heavy-tailed claim sizes with phase type ones and $\repsup[,i]\stackrel{d}{=}\repsup$, then a series expansion of the exact ruin probability is
  \begin{align*}
   \ruinmix     = \ruinrep
                   &+ \theta \sum_{n=1}^\infty \left(\frac{\epsilon}{1-\delta}\right)^n
                        \sum_{k=0}^{n-1} {n-1 \choose k} \theta^{k}(-\delta)^{n-1-k} \big( \disf{L_{n,k+1,n-1-k}} - \disf{L_{n-1,k,n-1-k}}\big) \notag\\
                   &- \delta \sum_{n=1}^\infty \left(\frac{\epsilon}{1-\delta}\right)^n
                        \sum_{k=0}^{n-1} {n-1 \choose k} \theta^{k}(-\delta)^{n-1-k} \big( \disf{L_{n,k,n-k}} - \disf{L_{n-1,k,n-1-k}}\big),
  \end{align*}
  where $\disf{L_{s,m,r}} = \pr(\repsup[,0]+\repsup[,1]+\dots+\repsup[,s]+\ehtc[1]+\dots+\ehtc[m] + \eptc[1]+\dots+\eptc[r] > u)$ and $\disf{L_{0,0,0}}=\ruinrep$. A sufficient condition for the convergence of the replace series expansion for all values of $u$ is $\epsilon<\lvert 1-\delta \rvert /\max\{\delta,\theta\}$.
\end{theorem}

Note that Theorem~\ref{t.replace expansion} gives only a sufficient condition for the convergence of the replace series expansion. If all parameters involved are explicitly known, one can find a necessary condition in the way indicated in the proof of Theorem~\ref{t.replace expansion}. In the next section, we propose two explicit approximations for the ruin probability based on these series expansions.

\section{Corrected phase-type approximations of the ruin probability}\label{S.Approximations}
The goal of this section is to provide approximations that maintain the numerical tractability but improve the accuracy of the phase-type approximations and that are able to capture the tail behavior of the exact ruin probability. Large deviations theory suggests that a single catastrophic event, i.e.\ a heavy-tailed stationary claim size \ehtc, is sufficient to cause ruin \cite{embrechts-MEE}. Observe that, for both the discard and replace series expansions, the second term contains a single appearance of \ehtc. For this reason, the proposed approximations for the ruin probability are constructed by the first two terms of their respective series expansions for the ruin probability (see Theorems~\ref{t.discard expansion} and \ref{t.replace expansion}), where the second term of each approximation is referred to as its correction term. We have the following definitions for the proposed approximations.

\begin{definition}\label{d.corrected discard}
  The corrected discard approximation of exact ruin probability \ruinmix\ is defined as
  \begin{equation}\label{e.corrected discard}
    \cordis := \ruindis + \frac{\epsilon\theta}{1-\delta+\epsilon\delta}\big(\pr(\dissup[,0]+\dissup[,1]+\ehtc[1] > u) - \pr(\dissup[,0] > u)\big),
  \end{equation}
  where \ruindis\ is the discard phase-type approximation of \ruinmix.
\end{definition}

In a similar manner, we define the approximation that connects to the replace expansion.

\begin{definition}\label{d.corrected replace}
  The corrected replace approximation of the exact ruin probability \ruinmix\ is given by the formula
  \begin{align}
    \correp := \ruinrep &+ \frac{\epsilon\theta}{1-\delta}\big(\pr(\repsup[,0]+\repsup[,1]+\ehtc[1]> u) - \pr(\repsup[,0] > u)\big)\notag \\
                               &- \frac{\epsilon\delta}{1-\delta}\big(\pr(\repsup[,0]+\repsup[,1]+\eptc[1]> u) - \pr(\repsup[,0] > u)\big), \label{e.corected replace}
  \end{align}
  where \ruinrep\ is the replace phase-type approximation of \ruinmix.
\end{definition}

In the following sections, we study characteristics of the corrected discard and the corrected replace approximations.

\subsection{Approximation errors}\label{SS.Approximation errors}
Due to the construction of the two corrected phase-type approximations, the discard and the replace, their difference from the exact ruin probability is the sum of the remaining terms, namely the terms for $n\geq 2$. For the error of the corrected discard approximation, we have the following theorem.

\begin{theorem}\label{t.bounds discard}
  The error of the corrected discard approximation is bounded from above and below as follows:
  \begin{equation*}
    \left(\frac{\epsilon\theta}{1-\delta+\epsilon\delta}\right)^2 \big(\disf{L_{\epsilon,2}} - \disf{L_{\epsilon,1}} \big) \quad
    \leq \ruinmix - \cordis \quad
    \leq \left(\frac{\epsilon\theta}{1-\delta+\epsilon\delta}\right)^2.
  \end{equation*}
\end{theorem}

\begin{remark}\label{r.error discard}
  Theorem~\ref{t.bounds discard} shows that the corrected discard approximation always underestimates the exact ruin probability, and its error is $O(\epsilon^2)$. Thus, the corrected discard approximation is a lower bound for the exact ruin probability.
\end{remark}

As done in the proof of Theorem~\ref{t.bounds discard}, similar probabilistic interpretations can also be given to the terms of the replace series expansion. However, due to the sign changes in the formula of the replace expansion (see Theorem~\ref{t.replace expansion}), it is not immediate whether the corrected replace approximation underestimates or overestimates the exact ruin probability. This depends on the characteristics of the distributions involved. As we see in Section~\ref{S.Example}, both overestimation and underestimation are possible. Studying the areas of over- or underestimation of the ruin probability is beyond the scope of this paper. In the sequel, we provide only absolute error bounds for the corrected replace approximation.

There are many possible ways to bound the error of the corrected replace approximation. For example, one could ignore all negative terms for $n\geq 2$ in the replace expansion and bound all positive terms. Of course, different techniques give different bounds. Among the different bounds we found, we present in Theorem~\ref{t.upper bound replace} the one that is valid for the biggest range of the perturbation parameter $\epsilon$.

\begin{theorem}\label{t.upper bound replace}
  When $\epsilon< \lvert 1-\delta\rvert/(\delta+\theta)$, an upper bound for the absolute error that we achieve with the corrected replace approximation is
  \begin{equation*}
    \lvert \ruinmix - \correp \rvert \leq \left( \frac{\epsilon}{1-\delta}\right)^2 (\delta+\theta)^2 \frac{1-\delta}{1-\delta-\epsilon(\delta+\theta)}.
  \end{equation*}
\end{theorem}

\begin{remark}
  Theorem~\ref{t.upper bound replace} shows that the absolute error of the replace approximation is $O(\epsilon^2)$. Note that the expression
  \begin{align*}
     \left(\frac{\epsilon}{1-\delta}\right)^2 \sum_{k=0}^{1} \theta^{k}(-\delta)^{1-k}
                \left[\theta \big( \disf{L_{2,k+1,1-k}} - \disf{L_{1,k,1-k}} \big)
                    - \delta \big( \disf{L_{2,k,2-k}} - \disf{L_{1,k,1-k}} \big)\right],
  \end{align*}
  which corresponds to the term of the replace expansion (see Theorem~\ref{t.replace expansion}) for $n=2$, is $O(\epsilon^2)$ and it could be used alternatively as an approximation of the real error.
\end{remark}

An advantage of the corrected discard approximation over the corrected replace is the following. The fact that the corrected discard approximation underestimates the exact ruin probability gives a positive sign for its error, namely its difference from the exact ruin probability, which according to Theorem~\ref{t.bounds discard} is bounded from above and below. This information with respect to the nature of its error makes the corrected discard approximation much more controllable than the corrected replace approximation. In the next section, we study the tail behavior of both corrected phase-type approximations.

\subsection{Tail behavior}\label{SS.Tail behavior}
To study the tail behavior of the two approximations, we assume that the distribution of $C^e$ belongs to the class of subexponential distributions $\mathcal{S}$. Following \cite{teugels75}, we give the following definition of $\mathcal{S}$.

\begin{definition}\label{d.subexponentiality}
   A distribution $F$ concentrated on $[0,\infty)$ belongs to the class of subexponential distributions $\mathcal{S}$ if and only if
   \begin{equation*}
     \lim_{u \rightarrow \infty} \frac{1-\disf{\con{F}{}}}{1-\disf{F}}=n, \qquad n=1,2,\dots
   \end{equation*}
\end{definition}

We use the notation $\df{f}\sim \df{g}$ to describe the relation $\lim_{u\rightarrow \infty}\df{f}/\df{g}=1$. When a distribution $F$ belongs to $\mathcal{S}$, it is known that $\com{F}$ decays slower than any exponential rate \cite{asmussen-RP}. Two very useful known properties of subexponentiality are the following, which are given without proof (see \cite{asmussen-RP}).

\begin{property}\label{p.subexponential class closed}
  The class $\mathcal{S}$ is closed under tail-equivalence. That is, if $\disf{\com{A}} \sim a \disf{\com{F}}$ for some $F \in \mathcal{S}$ and some constant $a > 0$, then $A \in \mathcal{S}$.
\end{property}

\begin{property}\label{p.lighter tail subexponential}
  Let $F \in \mathcal{S}$ and let $A$ be any distribution with a lighter tail, i.e. $\disf{\com{A}} = o\left(\disf{\com{F}}\right)$. Then for the convolution $A * F$ of $A$ and $F$ we have $A * F \in \mathcal{S}$ and $\disf{\com{\left(A * F\right)}} \sim \disf{\com{F}}$.
\end{property}

Before studying the tail behavior of the approximations, we first give the tail behavior of the exact ruin probability in the next theorem. We use the convention $A \in \mathcal{S}$ if the distribution of the r.v.\ $A$ belongs to $\mathcal{S}$.

\begin{theorem}\label{t.tail perturbed ruin probability}
  When $\ehtc \in \mathcal{S}$, the exact ruin probability \ruinmix\ has the following tail behavior:
  \begin{equation*}
    \ruinmix \sim \frac{\epsilon\theta}{1-\delta +\epsilon\delta-\epsilon\theta}\pr[(\ehtc > u)].
  \end{equation*}
\end{theorem}

For the tail behavior of the corrected discard approximation, the following result holds.
\begin{theorem}\label{t.tail discard}
  When $\ehtc \in \mathcal{S}$, we have for the corrected discard approximation the following tail behavior:
  \begin{equation}
    \cordis \sim \frac{\epsilon\theta}{1-\delta+\epsilon\delta} \pr[(\ehtc > u)].
  \end{equation}
\end{theorem}

Theorem~\ref{t.tail discard} shows that the corrected discard approximation captures the heavy-tailed behavior of the exact ruin probability, but is off by a term $\epsilon\theta$ in the denominator. In fact, for all values of parameters, the tail of the discard approximation is always below the tail of the exact ruin probability, which is expected since the discard approximation gives an underestimation of the exact ruin probability.

On the other hand, for the tail behavior of the corrected replace approximation, the following result holds.
\begin{theorem}\label{t.tail replace}
  When $\ehtc \in \mathcal{S}$, we have for the corrected replace approximation the following tail behavior:
  \begin{equation}
    \correp \sim \frac{\epsilon\theta}{1-\delta} \pr[(\ehtc > u)].
  \end{equation}
\end{theorem}

Comparing the coefficients of $\pr[(\ehtc > u)]$ in Theorems~\ref{t.tail discard} and \ref{t.tail replace}, we observe that the tail of the corrected replace approximation is always above the tail of the corrected discard approximation. To compare the tail behavior of the corrected replace approximation to that of the exact ruin probability, we only need to compare the coefficients of $\pr[(\ehtc > u)]$, and more precisely their denominators, as the expression with the largest denominator converges to zero faster. Therefore, the tails have the same behavior when $\e{\ptc}=\e{\htc}$, while the tail of the corrected replace approximation is above the tail of the exact ruin probability when $\e{\ptc}>\e{\htc}$ and below when $\e{\ptc}<\e{\htc}$.

\subsection{Relative error}\label{SS.Relative error}
Following the results of Section~\ref{SS.Tail behavior}, we show that the relative error at the tail for both approximations is $O(\epsilon)$.

\begin{lemma}\label{l.discard relative tail}
  When $\ehtc \in \mathcal{S}$, the relative error at the tail of the corrected discard approximation is
  \begin{equation*}
    \disf{R_{d,\epsilon}} = 1 - \frac{\cordis}{\ruinmix} \rightarrow \frac{\epsilon\theta}{1-\delta+\epsilon\delta}, \qquad \text{as} \quad u\rightarrow \infty.
  \end{equation*}
\end{lemma}

Recall that for the corrected replace approximation, different values of parameters lead to both over- and underestimation of the exact ruin probability. Thus, for this approximation it is more appropriate to evaluate the absolute relative error at its tail.

\begin{lemma}\label{l.replace relative tail}
  When $\ehtc \in \mathcal{S}$, the absolute relative error at the tail of the corrected replace approximation is
  \begin{equation}
    \left| \disf{R_{r,\epsilon}} \right|  = \left|  1 - \frac{\correp}{\ruinmix} \right|  \rightarrow \left| \frac{\epsilon(\theta-\delta)}{1-\delta} \right|, \qquad \text{as} \quad u\rightarrow \infty,
  \end{equation}
  and it goes asymptotically to zero when $\e{\ptc}=\e{\htc}$.
\end{lemma}

\begin{remark}\label{r.fixed tail behavior}
  Lemmas~\ref{l.discard relative tail} and \ref{l.replace relative tail} indicate that the relative errors of both corrected phase-type approximations do not converge to 0 as $u \rightarrow \infty$. However, the approximations give the exact value of the ruin probability at the origin and have guaranteed bounds of the order $O(\epsilon^2)$ for all values of $u$. On the other hand, the asymptotic result of Theorem~\ref{t.tail perturbed ruin probability} has the correct tail behavior but it gives relatively inaccurate estimates of the ruin probability for small values of $u$ for some combinations of the involved parameters. In order to provide a compromise between our approximations and the asymptotic result of Theorem~\ref{t.tail perturbed ruin probability}, one can simply change the coefficients of the correction terms (see Definitions~\ref{d.corrected discard} and \ref{d.corrected replace}) to $\epsilon\theta/(1-\delta +\epsilon\delta-\epsilon\theta)$, so that their tail behavior matches the correct tail behavior. Of course, one should also multiply the first terms of the approximations with proper coefficients to obtain $\ruindis[0]=\ruinrep[0]=\ruinmix[0]$. Such adjustments will lead to approximations with relative error at the tail that is asymptotically equal to zero. Moreover, the approximations work well for all values of the involved parameters, but may give worse results for small values of $u$ when they are compared with the original corrected phase-type approximations.
\end{remark}

The fact that the discard approximation always underestimates the ruin probability raises the question if it is possible to develop a result for its relative error for arbitrary values of $u$. The next theorem, which can be seen as the main technical contribution of the paper, shows that this is indeed possible.

\begin{theorem}\label{t.discard relative error}
  When $\ehtc \in \mathcal{S}$, there exists an $\eta>0$, such that for all $\epsilon<\eta$, the relative error $\disf{R_{d,\epsilon}}$ of the discard approximation at the point $u$ can be bounded by
  \begin{equation*}
    \disf{R_{d,\epsilon}} \leq  \frac{\epsilon\theta}{1-\delta+\epsilon\delta} \disf{H_\epsilon}+ \epsilon^2 K,
  \end{equation*}
  with $\disf{H_\epsilon}=\left(\frac{P(\dissup[,0]+\dissup[,1]+\dissup[,2]+\ehtc[1]+\ehtc[2] > u)}{P(\dissup[,0]+\dissup[,1]+\ehtc[1] > u)}-1\right)$ and $K$ a finite constant.
\end{theorem}

The bound is sharp in the sense that $\disf{H_\epsilon}\rightarrow 1$ as $u\rightarrow\infty$, which recovers the relative error at the tail, up to a term $O(\epsilon^2)$. Moreover, $\disf{H_\epsilon}$ is uniformly bounded in $u$ and $\epsilon$.

\section{Numerical examples}\label{S.Example}
In Section~\ref{S.Power series expansion}, we pointed out that the first terms of the discard and the replace expansions are phase-type approximations of \ruinmix. The goal of this section is to show numerically that adding the second term of these expansions leads to improved approximations (corrected discard and corrected replace approximations respectively) that are significantly more accurate than their phase-type counterparts. Moreover, the additional term has a great impact on the accuracy of the improved approximations even for small values of the perturbation parameter.

Therefore, in this section we check the accuracy of the corrected discard (see Definition~\ref{d.corrected discard}) and the corrected replace approximations (see Definition~\ref{d.corrected replace}) by comparing them with the exact ruin probability and their corresponding phase-type approximations. Since it is more meaningful to compare approximations with exact results than with simulation outcomes, we choose the general claim size distributions $G_\epsilon$ such that there exists an exact formula for the ruin probability \ruinmix.

In Section~\ref{SS.Theorem exact perturebed ruin probability}, we derive the exact formula for the ruin probability \ruinmix\ for a specific choice of the claim size distribution. Using the latter claim size distribution, in Section~\ref{SS.Numerical results} we perform our numerical experiments and we draw our conclusions.

\subsection{Test distribution}\label{SS.Theorem exact perturebed ruin probability}
As claim size distribution we use a mixture of an exponential distribution with rate $\nu$ and a heavy-tailed one that belongs to a class of long-tailed distributions introduced in \cite{abate99a}. The Laplace transform of the latter distribution is $\lthtc = 1-\frac{s}{(\mu+\sqrt{s})(1+\sqrt{s})}$, where $\e C=\mu^{-1}$ and all higher moments are infinite. Furthermore, the Laplace transform of the stationary heavy-tailed claim size distribution is
\begin{equation*}
  \ltehtc = \frac{\mu}{(\mu+\sqrt{s})(1+\sqrt{s})},
\end{equation*}
which for $\mu\neq 1$ can take the form
\begin{equation*}
  \ltehtc = \left(\frac{\mu}{1-\mu}\right) \left(\frac{1}{\mu+\sqrt{s}} - \frac{1}{1+\sqrt{s}}\right).
\end{equation*}

For this combination of claim size distributions, the ruin probability can be found explicitly:

\begin{theorem}\label{t.perturbed ruin probability example}
  Assume that claims arrive according to a Poisson process with rate $\lambda$, the premium rate is 1 and the Laplace transform of the claim size distribution is
  \begin{equation}\label{e.perturbed claim size distribution example}
    \ltgc = (1-\epsilon) \frac{\nu}{s+\nu} + \epsilon\left(1-\frac{s}{(\mu+\sqrt{s})(1+\sqrt{s})}\right),
  \end{equation}
  with $\rho_\epsilon = \frac{\lambda}{\mu\nu} \big(\mu + \epsilon(\nu-\mu)\big)<1$. For this mixture model, the ruin probability is
  \begin{equation}\label{e.perturbed ruin probability example}
    \ruinmix = \frac{\lambda}{\mu\nu} \Big(\mu\nu - \lambda \big(\mu + \epsilon(\nu-\mu)\big) \Big)\sum_{i=1}^4 \frac{a_i}{\nu_i(\epsilon)} \df[\nu_i^2(\epsilon) u]{\zeta},
  \end{equation}
  where
  \begin{equation}\label{e.error function}
    \df{\zeta} := e^u\erfct{\sqrt{u}},
  \end{equation}
  and $-\nu_i(\epsilon)$, $i=1,\dots,4$, are the roots of the polynomial
  \begin{equation*}
    \df[x]{d} = x^4 +(\mu+1) x^3 + (\mu+\nu-\lambda) x^2 + (\mu+1)(\nu - \lambda+\lambda\epsilon)x + \big(\mu(\nu-\lambda)+\lambda\epsilon(\mu-\nu)\big).
  \end{equation*}
  Finally, the coefficients $a_i$ satisfy $a_i=\lim_{x\rightarrow -\nu_i(\epsilon)} \frac{\df[x]{n}}{\df[x]{d}}\big(x+\nu_i(\epsilon)\big)$, $i=1,\dots,4$, where
  \begin{equation*}
    \df[x]{n} = (1-\epsilon)(\mu+x)(1+x) +  \epsilon (x^2+\nu).
  \end{equation*}
\end{theorem}

\subsection{Numerical results}\label{SS.Numerical results}
In this section, we fix values for the parameters of the mixture model described in the previous section and we perform our numerical experiments. Although we do not have any restrictions for the parameters of the involved claim size distributions, from a modeling point of view, it is counterintuitive to fit a heavy-tailed claim size distribution with a mean smaller than the mean of the phase-type claim size distribution. For this reason, we select $\mu=2$ and $\nu=3$.

For the perturbation parameter $\epsilon$, the only restrictions arise from the conditions for the convergence of the discard and the replace series expansions (see Appendix) and the stability condition. A closer look at the formulas reveals that, in the case of unequal means, for every value of $\epsilon$ there exists a value for the arrival rate $\lambda$ such that all conditions are satisfied. However, a logical constraint for the perturbation parameter is $\epsilon\leq 0.1$. The reason for this constraint is that in the case of phase-type approximations it is not natural to remove more than 10\% of the data.

To start our experiments, we first choose the ``worst case scenario" for the perturbation parameter, which is $\epsilon=0.1$. It seems that this ``worst case scenario" for the perturbation parameter is the ``best case scenario" for the improvement we can achieve with the corrected phase-type approximations. When the perturbation parameter is big enough, a lot of information with respect to the tail behavior of the ruin probability is missing from its phase-type approximations. So, it is quite natural to expect a great improvement when we add the second term of the respective series expansion, which contains a big part of this missing information. In this scenario, we compare the corrected phase-type approximations with their respective phase-type approximations when $\rho_{0.1}$ takes the values 0.5, 0.7 and 0.9.

\begin{figure}
         \centering
         \begin{subfigure}[]
                 \centering
                 \includegraphics[scale=1.0]{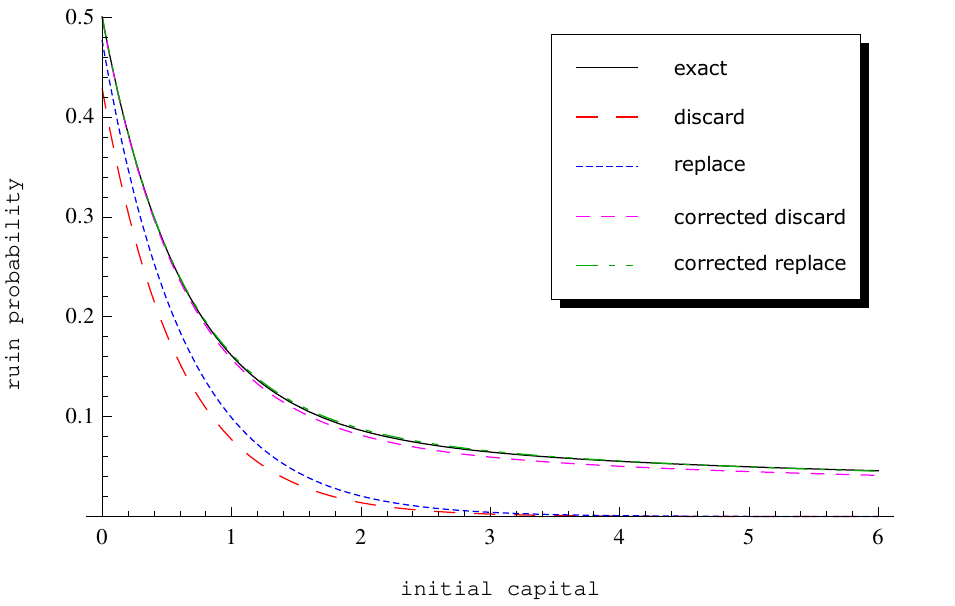}
                 \label{fig:A}
         \end{subfigure}%
         \begin{subfigure}[]
                 \centering
                 \includegraphics[scale=1.0]{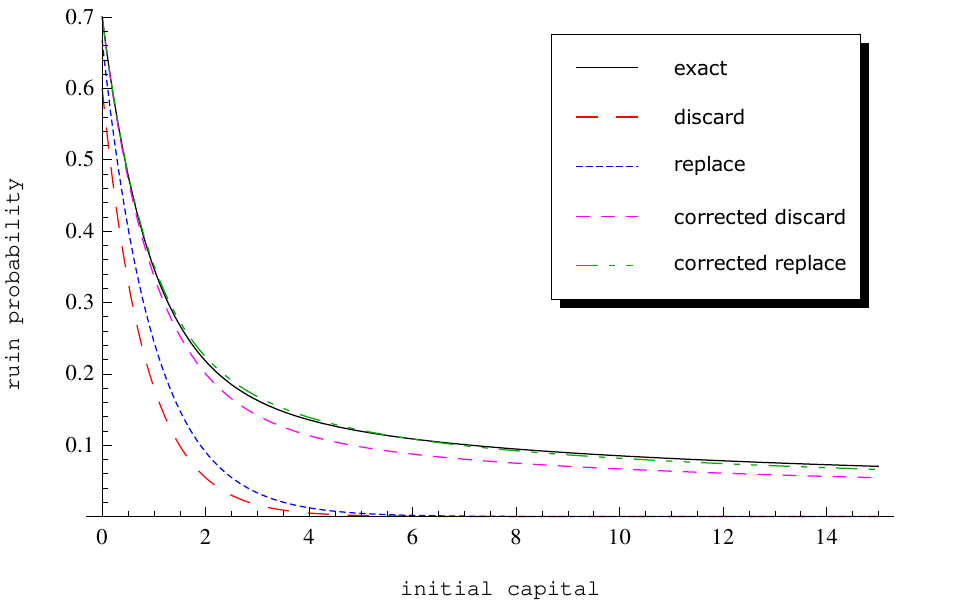}
                 \label{fig:B}
         \end{subfigure}
         \begin{subfigure}[]
                 \centering
                 \includegraphics[scale=1.0]{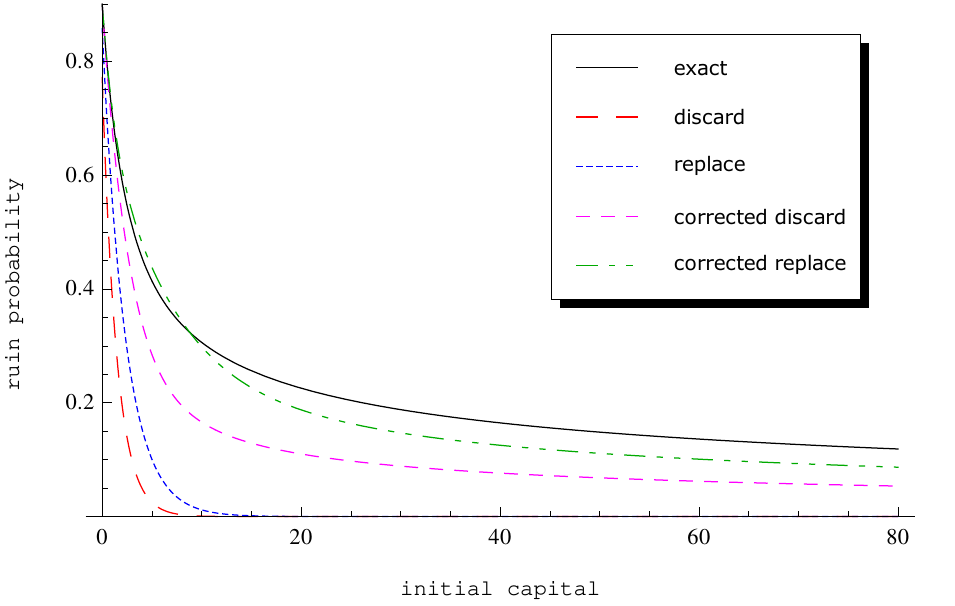}
                 \label{fig:C}
         \end{subfigure}
         \caption{Exact ruin probability with phase-type and corrected phase-type approximations for perturbation parameter 0.1 and average claim rate: (a) 0.5, (b) 0.7, and (c) 0.9.}\label{f.approximations}
 \end{figure}

From Figure~\ref{f.approximations}, we conclude that the corrected discard and the corrected replace approximations manage to reduce the ``gap" between their respective phase-type approximations and the exact ruin probability. Although the scale of the graphs is different, it is evident that the gap closes more efficiently for small values of $\rho_\epsilon$, a conclusion that can be also supported theoretically by Section~\ref{SS.Tail behavior}. Furthermore, the corrected replace approximation overestimates the ruin probability for small values of $u$ and, as expected, it is better at the tail than the corrected discard approximation.

For small values of $\rho_\epsilon$ and small values of $\epsilon$, one could argue that the gap between the exact ruin probability and its phase-type approximations is so small that the corrected phase-type approximations do not improve on the accuracy of their phase-type counterparts. For this reason, we choose $\epsilon=0.001$ and $\rho_{0.001}=0.5$, and we compare all approximations with the exact ruin probability. We show that the improvement we achieve with the corrected phase-type approximations is still significant, even for this seemingly ``bad scenario".

\begin{table}[h]
  \begin{center}
    \begin{tabular}{|c|ccccc|}\hline
      $u$   &exact ruin probability     & discard                     & replace                      & corrected discard            & corrected replace \\ \hline
      0     & 0.50000000                & 0.49925037                  & 0.49975012                   & 0.50000000                   & 0.50000000  \\
      1     & 0.11211000                & 0.11114757                  & 0.11142576                   & 0.11210955                   & 0.11211017  \\
      2     & 0.02557910                & 0.02474466                  & 0.02484381                   & 0.02557847                   & 0.02557930  \\
      3     & 0.00621454                & 0.00550887                  & 0.00553925                   & 0.00621386                   & 0.00621466 \\
      4     & 0.00184042                & 0.00122643                  & 0.00123504                   & 0.00183975                   & 0.00184047 \\
      5     & 0.00082276                & 0.00027304                  & 0.00027536                   & 0.00082212                   & 0.00082275  \\
      6     & 0.00056334                & 0.00006078                  & 0.00006139                   & 0.00056273                   & 0.00056329  \\
      7     & 0.00047969                & 0.00001353                  & 0.00001368                   & 0.00047910                   & 0.00047962  \\
      8     & 0.00043993                & $3.01 \times 10^{-6}$       & $3.05 \times 10^{-6}$        & 0.00043937                   & 0.00043985  \\
      9     & 0.00041336                & $6.70 \times 10^{-7}$       & $6.80 \times 10^{-7}$        & 0.00041284                   & 0.00041329  \\
      10    & 0.00039235                & $1.49 \times 10^{-7}$       & $1.51 \times 10^{-7}$        & 0.00039183                   & 0.00039225  \\ \hline
    \end{tabular}
  \end{center}
  \caption{Exact ruin probability with phase-type and corrected phase-type approximations for perturbation parameter 0.001 and average claim rate 0.5.}\label{table-error}
\end{table}

From Table~\ref{table-error}, we observe that even for this small value of $\epsilon$ the corrected discard and the corrected replace approximations yield significant improvements for their respective phase-type approximations. The difference between the exact ruin probability and the corrected phase-type approximations is $O(10^{-6})$, while for the phase-type approximations it is $O(10^{-3})$. In order to understand the magnitude of the improvement we achieve with the corrected phase-type approximations we need to look also at the relative errors of all the approximations involved. It is evident that the relative error of the phase-type approximations easily reaches values close to 1 (approximately after value 5 of the initial capital in this example), while the corrected phase-type approximations give a relative error $O(\epsilon)$.

An interesting observation is that the corrected replace approximation gives better numerical estimations than the corrected discard approximation, both in absolute and relative errors. However, due to the sign changes in the formula of the replace expansion (see Theorem~\ref{t.replace expansion}) it is difficult to find tight bounds for this approximation.

Finally, note that we performed extensive numerical experiments for various values of the perturbation parameter $\epsilon$ in the interval $[0.001,0.1]$. We chose to present only the extreme cases, since the qualitative conclusions for the intermediate values of $\epsilon$ are similar to those of the extreme cases.

\section{Total loss and Value at Risk}\label{S.VaR}
In this section, we give a brief overview of how our technique works when we calculate quantities in finite time horizon. As test example, we use the aggregate loss in a fixed period, and we provide the corrected phase-type approximations when the aggregate loss is a compound Poisson sum. Moreover, we extend our technique in case the aggregate loss is a compound mixed Poisson sum. Finally, we perform a small numerical experiment to compare the Value at Risk (VaR) for a given level $\alpha$ that we obtain from the original distribution, the corrected phase-type approximation and its corresponding phase-type approximation.

Suppose that we are interested in evaluating the aggregate loss in a fixed period $[0,t]$. The number \df[t]{N_\epsilon} of claims \mixc over this fixed period follows a Poisson distribution with rate $\lambda t$. Observe that $\df[t]{N_\epsilon}$ can be seen as a superposition of two independent Poisson processes \df[t]{N_\epsilon^P} and \df[t]{N_\epsilon^H}, with rates $\lambda(1-\epsilon)t$ and $\lambda \epsilon t$ for the phase-type and the heavy-tailed claims sizes, respectively. Thus, we write
\begin{equation*}
  \loss := \pr\Bigg(\sum_{i=1}^{\df[t]{N_\epsilon}} \mixc[,i] >x \Bigg) = \pr\Bigg(\sum_{i=1}^{\df[t]{N_\epsilon^P}} \ptc[i] + \sum_{i=1}^{\df[t]{N_\epsilon^H}} \htc[i] >x \Bigg).
\end{equation*}

To find \loss, we condition on the number \df[t]{N_\epsilon^H} of the heavy-tailed claim sizes and we get
\begin{align*}
  \loss =&\pr\Bigg(\sum_{i=1}^{\df[t]{N_\epsilon^P}} \ptc[i] + \sum_{i=1}^{\df[t]{N_\epsilon^H}} \htc[i] >x \Bigg)
  = \sum_{k=0}^\infty \pr\Bigg(\sum_{i=1}^{\df[t]{N_\epsilon^P}} \ptc[i] + \sum_{i=1}^k \htc[i] >x \Bigg) \pr \big(\df[t]{N_\epsilon^H} = k \big) \\
  =& \pr\Bigg(\sum_{i=1}^{\df[t]{N_\epsilon^P}} \ptc[i] >x \Bigg) \pr \big(\df[t]{N_\epsilon^H} = 0 \big)
   +\pr\Bigg(\sum_{i=1}^{\df[t]{N_\epsilon^P}} \ptc[i] + \htc >x \Bigg) \pr \big(\df[t]{N_\epsilon^H} = 1 \big) \\
   &+\pr\Bigg( \left. \sum_{i=1}^{\df[t]{N_\epsilon^P}} \ptc[i] + \sum_{i=1}^{\df[t]{N_\epsilon^H}} \htc[i]  >x \right|  \df[t]{N_\epsilon^H} \geq 2 \Bigg) \pr \big(\df[t]{N_\epsilon^H} \geq 2 \big)=\\
  =& \underbrace{\pr\Bigg(\sum_{i=1}^{\df[t]{N_\epsilon^P}} \ptc[i] >x \Bigg)}_{\text{PH-approximation}} e^{-\lambda \epsilon t}
   + \pr \big(\df[t]{N_\epsilon^H} \geq 1 \big) \pr\Bigg(\sum_{i=1}^{\df[t]{N_\epsilon^P}} \ptc[i] + \htc >x \Bigg)  \\
   &+ \pr \big(\df[t]{N_\epsilon^H} \geq 2 \big) \Bigg[\pr\Bigg( \left. \sum_{i=1}^{\df[t]{N_\epsilon^P}} \ptc[i] + \sum_{i=1}^{\df[t]{N_\epsilon^H}} \htc[i]  >x \right|  \df[t]{N_\epsilon^H} \geq 2 \Bigg) - \pr\Bigg(\sum_{i=1}^{\df[t]{N_\epsilon^P}} \ptc[i] + \htc >x \Bigg) \Bigg].
\end{align*}

As in the case of ruin probabilities, we define the corrected phase-type approximation by keeping only the terms that contain at most one appearance of the heavy-tailed claim sizes. Thus, we have the following definition.
\begin{definition}\label{d.loss corrected discard approximation}
  The corrected discard approximation of the tail of the aggregated claim sizes in a fixed time interval $[0,t]$ is defined as
  \begin{equation*}
    \cordisloss := e^{-\lambda \epsilon t} \pr\Bigg(\sum_{i=1}^{\df[t]{N_\epsilon^P}} \ptc[i] >x \Bigg)
   + (1-e^{-\lambda \epsilon t}) \pr\Bigg(\sum_{i=1}^{\df[t]{N_\epsilon^P}} \ptc[i] + \htc >x \Bigg) .
  \end{equation*}
\end{definition}

Observe that the coefficient of the correction term in Definition~\ref{d.loss corrected discard approximation} is equal to $\pr \big(\df[t]{N_\epsilon^H} \geq 1 \big)$ and not $\pr \big(\df[t]{N_\epsilon^H} = 1 \big)$ as one would expect. We used this modification in order to achieve more accurate estimates of the aggregate loss, without loosing the main characteristic of the corrected discard approximation, which is the fact that it underestimates the exact distribution. According to the next theorem, the approximation error of \cordisloss\ is of order $O(\epsilon^2)$. As it was the case for Theorem~\ref{t.upper bound replace}, there are many ways to find a lower bound for the error. In the next theorem, we present a bound yielding a simple expression.
\begin{theorem}\label{t.loss corrected discard error bounds}
  The error of the corrected discard approximation \cordisloss\ is bounded as follows:
  \begin{equation*}
    \pr \big(\df[t]{N_\epsilon^H} \geq 2 \big) \Bigg[\pr\Bigg( \sum_{i=1}^{\df[t]{N_\epsilon^P}} \ptc[i] + \htc[1] + \htc[2]  >x \Bigg) - \pr\Bigg(\sum_{i=1}^{\df[t]{N_\epsilon^P}} \ptc[i] + \htc[1] >x \Bigg) \Bigg] \leq \loss - \cordisloss  \leq  \epsilon^2 (\lambda t)^2.
  \end{equation*}
\end{theorem}

\begin{remark}\label{r.loss corrected replace}
  The corrected replace approximation can be constructed in a similar manner. However, special attention should be paid to the fact that we need to condition not only on the number \df[t]{N_\epsilon^H} of heavy-tailed claims but also on the total number of claims, namely \df[t]{N_\epsilon}. This of course will lead to expressions with the same order of complexity with that of the approximation in Definition~\ref{d.corrected replace}.
\end{remark}

If the time $t$ we are interested in is not fixed but a random variable, e.g.\ $T$, the total aggregate loss is a compound mixed Poisson r.v. The corrected discard approximation takes the form
\begin{equation*}
  \cordisloss[T] = \int_{0}^\infty e^{-\lambda \epsilon t} \pr\Bigg(\sum_{i=1}^{\df[t]{N_\epsilon^P}} \ptc[i] >x  \Bigg) d\pr(T\leq t)
                    + \int_{0}^\infty (1-e^{-\lambda \epsilon t}) \pr\Bigg(\sum_{i=1}^{\df[t]{N_\epsilon^P}} \ptc[i] + \htc >x \Bigg) d\pr(T\leq t),
\end{equation*}
and an upper bound for its error is $\epsilon^2 \lambda^2 \e T^2$. As a last result, we find a compact formula for the Laplace-Stieltjes transform of the \cordisloss[T]. We use the notation \ltptc{}, \lthtc\ and \lttime for the Laplace transforms of the phase-type claim sizes, the heavy-tailed claim sizes and the r.v.\ $T$, respectively.

\begin{theorem}\label{t.LST mixed poisson}
   The Laplace-Stieltjes transform of \cordisloss[T] is given by the formula
  \begin{equation*}
    \mathcal{L}\{\cordisloss[T]\} = \frac1s -\frac{1 - \lthtc}s \lttime[\lambda \big(1- (1-\epsilon) \ltptc{}\big)] - \frac{\lthtc}s \lttime[\lambda (1-\epsilon) \big(1-\ltptc{}\big)].
  \end{equation*}
\end{theorem}

One can find the corrected discard approximation analytically (or numerically) by applying Laplace inversion to $\mathcal{L}\{\cordisloss[T]\}$.

A widely used risk measure that connects to the aggregate loss, is the Value at Risk (VaR), which is defined as the threshold value such that the probability of the aggregate loss to exceed this value is less than a given level $\alpha$. In other words, the VaR is equal to the $(1-\alpha)$-quantile of \loss. We show through a small numerical experiment that the VaR that is estimated with the corrected discard approximation is closer to the original VaR, than the one we obtain with the discard phase-type approximation. For our example, we choose the arrival rate $\lambda=1$, the service time distribution a mixture of an exponential distribution with rate $3/2$ and a Pareto distribution with scale and shape parameters 1 and 2 respectively, and $\epsilon=0.01$. We estimate the VaR values at level $0.99$ for the interval $[0,t]$, for the values of $t=1,5,10,15,20$. Note that we simulated the system in order to estimate the exact VaR values. We summarize our results in Table~\ref{table-var}.

\begin{table}[h]
  \begin{center}
    \begin{tabular}{|c|ccc|}\hline
       $t$  & simulation    & discard          & corrected discard  \\ \hline
        1   & 4.16          & 4.09              & 4.14 \\
        5   & 9.77          & 9.54              & 9.74 \\
        10  & 15.24         & 14.89             & 15.22 \\
        15  & 20.27         & 19.73             & 20.17 \\
        20  & 24.99         & 23.76             & 24.30 \\ \hline
    \end{tabular}
  \end{center}
  \caption{Comparison of VaR values we obtain from simulation results and the phase-type and corrected phase-type approximations. The VaR level is $\alpha=0.99$, the perturbation parameter 0.01 and the average claim rate 0.67.}\label{table-var}
\end{table}

We want to point out here, that this numerical study differs from our previous examples. Although in all other examples we were comparing tail probabilities at given values, here we compare the values at which the original distribution and its approximations give us the same tail probability. This observation explains why the difference between the values in Table~\ref{table-var} are not of order $O(\epsilon^2)$.

\section*{Acknowledgments}\label{acknowledgments}
The work of Maria Vlasiou and Eleni Vatamidou is supported by Netherlands Organisation for Scientific Research (NWO) through project number 613.001.006. The work of Bert Zwart is supported by an NWO VIDI grant and an IBM faculty award. Bert Zwart is also affiliated with VU University Amsterdam, and the Georgia Institute of Technology. The authors thank the referee for a careful reading and several constructive comments, of which one led to the example in Section~\ref{S.VaR}.

\section*{Appendix}\label{appendix}

\begin{proof}[Proof of Theorem~\ref{t.discard expansion}]
  From Eq.~\eqref{e.LT PK perturbed ruin probability} and \eqref{e.LT PK discard ruin probability} we find
    \begin{align*}
     \ltruinmix &= \frac{1-(1-\epsilon)\delta-\epsilon\theta}{1-(1-\epsilon)\delta \lteptc - \epsilon\theta \ltehtc}
                 = \frac{(1-(1-\epsilon)\delta)-\epsilon\theta}{\frac{1-(1-\epsilon)\delta}{\ltruindis} - \epsilon\theta \ltehtc}
                 = \frac{1-\frac{\epsilon\theta}{1-(1-\epsilon)\delta}}{\frac{1}{\ltruindis} - \frac{\epsilon\theta}{1-(1-\epsilon)\delta} \ltehtc}\\
                &= \ltruindis \left(1-\frac{\epsilon\theta}{1-\delta+\epsilon\delta}\right)\frac{1}{1-\frac{\epsilon\theta}{1-\delta+\epsilon\delta} \ltruindis \ltehtc}\\
                &= \ltruindis \left(1-\frac{\epsilon\theta}{1-\delta+\epsilon\delta}\right) \sum_{n=0}^\infty \left(\frac{\epsilon\theta}{1-\delta+\epsilon\delta}\right)^n
                   \left(\ltruindis \ltehtc\right)^n\\
                &= \ltruindis \sum_{n=0}^\infty \left(\frac{\epsilon\theta}{1-\delta+\epsilon\delta}\right)^n \left(\ltruindis \ltehtc\right)^n
                   - \sum_{n=0}^\infty \left(\frac{\epsilon\theta}{1-\delta+\epsilon\delta}\right)^{n+1} \left(\ltruindis\right)^{n+1} \left(\ltehtc\right)^n\\
                &= \ltruindis + \sum_{n=1}^\infty \left(\frac{\epsilon\theta}{1-\delta+\epsilon\delta}\right)^n
                   \left[\left(\ltruindis\right)^{n+1}\left(\ltehtc\right)^n - \left(\ltruindis\right)^n\left(\ltehtc\right)^{n-1}\right].
    \end{align*}
    Using Laplace inversion we obtain
    \begin{equation*}
      \ruinmix = \ruindis + \sum_{n=1}^\infty \left(\frac{\epsilon\theta}{1-\delta+\epsilon\delta}\right)^n \left(\disf{L_{\epsilon,n}} - \disf{L_{\epsilon,n-1}} \right),
    \end{equation*}
    where $\disf{L_{\epsilon,n}} = \pr(M^\bullet_{\epsilon,0}+M^\bullet_{\epsilon,1}+\dots+M^\bullet_{\epsilon,n}+C^e_1+\dots+C^e_n > u)$. Note that this power series expansion is valid if and only if
    \begin{equation*}
      \left| \frac{\epsilon\theta}{1-\delta+\epsilon\delta} \ltruindis \ltehtc \right| <1.
    \end{equation*}
    We know that $\lvert \ltruindis \ltehtc \rvert \leq 1$, so a necessary and sufficient condition for the convergence of the power series for all values of $s$ is $\lvert \epsilon\theta \rvert < \lvert 1-\delta+\epsilon\delta \rvert$. If we assume that $\theta > \delta$, then an immediate consequence of the stability condition $\rho_\epsilon<1$ is that $(1-\epsilon)\delta<1 $. Therefore the convergence condition simplifies to $\epsilon < (1-\delta)/(\theta -\delta)$.
\end{proof}

\begin{proof}[Proof of Theorem~\ref{t.replace expansion}]
  We set $\disf[s]{\lt D} = \theta \ltehtc -\delta \lteptc$. By using \eqref{e.LT PK perturbed ruin probability} and \eqref{e.LT PK replace ruin probability} we find
  \begin{align*}
        \ltruinmix &
                    = \frac{1-(1-\epsilon)\delta-\epsilon \theta}{1-(1-\epsilon)\delta \lteptc - \epsilon\theta \ltehtc}
                    = \frac{(1-\delta)-\epsilon(\theta-\delta)}{1-\delta \lteptc-\epsilon (\theta \ltehtc -\delta \lteptc)}
                    = \frac{(1-\delta)-\epsilon(\theta-\delta)}{\frac{1-\delta}{\ltruinrep} - \epsilon \disf[s]{\lt D} }\\
                   &= \frac{1-\epsilon\frac{\theta - \delta}{1-\delta}}{\frac{1}{\ltruinrep} - \frac{\epsilon}{1-\delta} \disf[s]{\lt D}}
                    = \ltruinrep \left(1-\epsilon\frac{\theta - \delta}{1-\delta}\right) \frac{1}{1-\frac{\epsilon}{1-\delta} \ltruinrep \disf[s]{\lt D}}\\
                   &= \ltruinrep \left(1-\epsilon\frac{\theta - \delta}{1-\delta}\right)
                      \sum_{n=0}^\infty \left(\frac{\epsilon}{1-\delta}\right)^n
                      \left(\ltruinrep \disf[s]{\lt D}\right)^n\\
                   &= \ltruinrep \sum_{n=0}^\infty \left(\frac{\epsilon}{1-\delta}\right)^n \left(\ltruinrep \disf[s]{\lt D}\right)^n
                      - (\theta-\delta)\sum_{n=0}^\infty \left(\frac{\epsilon}{1-\delta}\right)^{n+1} \left(\ltruinrep\right)^{n+1} \left(\disf[s]{\lt D}\right)^n\\
                   &= \ltruinrep + \sum_{n=1}^\infty \left(\frac{\epsilon}{1-\delta}\right)^n
                      \left(\ltruinrep\right)^n \left[\ltruinrep\left(\disf[s]{\lt D}\right)^n - (\theta-\delta)\left(\disf[s]{\lt D}\right)^{n-1}\right].
  \end{align*}
  But,
  \begin{align*}
      \ltruinrep \left(\disf[s]{\lt D}\right)^n &- (\theta-\delta)\left(\disf[s]{\lt D}\right)^{n-1}\\
                                 =& \ltruinrep \sum_{k=0}^n {n \choose k}(\theta \ltehtc)^k (-\delta \lteptc)^{n-k} - (\theta-\delta)\sum_{k=0}^{n-1} {n-1 \choose k}(\theta \ltehtc)^k (-\delta \lteptc)^{n-1-k}\\
                                 =& \ltruinrep \left[\sum_{k=1}^{n-1} \left({n-1 \choose k} + {n-1 \choose k-1}\right) (\theta \ltehtc)^k (-\delta \lteptc)^{n-k} + (\theta \ltehtc)^n + (-\delta \lteptc)^n\right]\\
                                  &- (\theta-\delta)\sum_{k=0}^{n-1} {n-1 \choose k}(\theta \ltehtc)^k (-\delta \lteptc)^{n-1-k}\\
                                 =& \theta \sum_{k=0}^{n-1} {n-1 \choose k}\left(\ltruinrep \ltehtc - 1\right) (\theta \ltehtc)^{k} (-\delta \lteptc)^{n-1-k}\\
                                  &- \delta \sum_{k=0}^{n-1} {n-1 \choose k}\left(\ltruinrep \lteptc - 1\right) (\theta \ltehtc)^{k} (-\delta \lteptc)^{n-1-k}.
  \end{align*}

  \noindent Therefore,
  \begin{align*}
      \ltruinmix =& \ltruinrep \\
                  &+ \theta \sum_{n=1}^\infty \left(\frac{\epsilon}{1-\delta}\right)^n \left(\left(\ltruinrep\right)^{n+1}\ltehtc - \left(\ltruinrep\right)^n\right)
                    \sum_{k=0}^{n-1} {n-1 \choose k} (\theta \ltehtc)^{k} (-\delta \lteptc)^{n-1-k}\\
                  &- \delta \sum_{n=1}^\infty \left(\frac{\epsilon}{1-\delta}\right)^n \left(\left(\ltruinrep\right)^{n+1} \lteptc - \left(\ltruinrep\right)^n\right)
                    \sum_{k=0}^{n-1} {n-1 \choose k} (\theta \ltehtc)^{k} (-\delta \lteptc)^{n-1-k}.
  \end{align*}
  Applying Laplace inversion we find
  \begin{align*}
     \ruinmix     = \ruinrep
                    &+ \theta \sum_{n=1}^\infty \left(\frac{\epsilon}{1-\delta}\right)^n
                        \sum_{k=0}^{n-1} {n-1 \choose k} \theta^{k}(-\delta)^{n-1-k} \left( \disf{L_{n,k+1,n-1-k}} - \disf{L_{n-1,k,n-1-k}}\right)\\
                    &- \delta \sum_{n=1}^\infty \left(\frac{\epsilon}{1-\delta}\right)^n
                        \sum_{k=0}^{n-1} {n-1 \choose k} \theta^{k}(-\delta)^{n-1-k} \left( \disf{L_{n,k,n-k}} - \disf{L_{n-1,k,n-1-k}}\right),
  \end{align*}
  where $\disf{L_{s,m,r}} = \pr(\repsup[,0]+\repsup[,1]+\dots+\repsup[,s]+\ehtc[1]+\dots+\ehtc[m] + \eptc[1]+\dots+\eptc[r] > u)$. Similarly to the discard expansion, the replace series converges for a given value of $s$ if and only if
  \begin{equation*}
    \left| \frac{\epsilon}{1-\delta}\ltruinrep\left(\theta \ltehtc -\delta \lteptc \right) \right| < 1.
  \end{equation*}
  If $\sigma = \max_s \lvert \ltruinrep\left( \theta \ltehtc -\delta \lteptc \right) \rvert$, then a necessary and sufficient condition for the convergence of the replace series for all values of $s$ is $\epsilon<\lvert 1-\delta \rvert /\sigma$. However, we do not have exact formulas for the Laplace transforms \ltruinrep, \ltehtc\ and \lteptc, and thus we can only find a sufficient condition for the convergence of the series. It can easily be shown that $\max_{s\geq 0} \lvert \theta \ltehtc -\delta \lteptc \rvert \leq \max\{\delta,\theta\}$. Since $\lvert \ltruinrep \rvert \leq 1$, a sufficient condition for the convergence of the replace series is $\epsilon<\lvert 1-\delta \rvert /\max\{\delta,\theta\}$. When $\theta>\delta$, the condition simplifies to $\epsilon< (1-\delta)/\theta$.
\end{proof}

\begin{proof}[Proof of Theorem~\ref{t.bounds discard}]
    An interesting observation is that we can interpret the terms $\disf{L_{\epsilon,n}}-\disf{L_{\epsilon,n-1}}$ in Theorem~\ref{t.discard expansion} in terms of a renewal process $\{N_{D,\epsilon}(u),u\geq 0\}$ with a delayed first renewal \dissup[,0]. Consequently, $\pr(N_{D,\epsilon}(u)=~0) =  \disf{L_{\epsilon,0}}$ and $\pr(N_{D,\epsilon}(u)= n) = \disf{L_{\epsilon,n}} - \disf{L_{\epsilon,n-1}}$, for $n\geq 1$. As a result,
    \begin{align*}
      \ruinmix - \cordis &= \sum_{n=2}^\infty \left(\frac{\epsilon\theta}{1-\delta+\epsilon\delta}\right)^n \pr(N_{D,\epsilon}(u) = n)\\
                                      &= \left(\frac{\epsilon\theta}{1-\delta+\epsilon\delta}\right)^2
                                       \e{\left[  \left(\frac{\epsilon\theta}{1-\delta+\epsilon\delta}\right)^{N_{D,\epsilon}(u)-2} \mathds{I}\left(N_{D,\epsilon}(u)\geq2\right)\right]}\\
                                     &\leq \left(\frac{\epsilon\theta}{1-\delta+\epsilon\delta}\right)^2,
    \end{align*}
    where the latter inequality holds because $\frac{\epsilon\theta}{1-\delta+\epsilon\delta}<1$. Thus, an upper bound for the approximation error is $\left(\frac{\epsilon\theta}{1-\delta+\epsilon\delta}\right)^2$. Due to the renewal argument, all terms in the discard series expansion are positive. Consequently, the corrected discard approximation always underestimates the exact ruin probability and the term $\left(\frac{\epsilon\theta}{1-\delta+\epsilon\delta}\right)^2 \big(\disf{L_{\epsilon,2}} - \disf{L_{\epsilon,1}} \big)$ is a lower bound for the achieved error.
\end{proof}

\begin{proof}[Proof of Theorem~\ref{t.upper bound replace}]
  Using the triangular inequality and the fact that the distance between two distributions is smaller than or equal to 1, we obtain
  \begin{align*}
    \lvert \ruinmix - \correp \rvert &\leq (\delta + \theta) \sum_{n=2}^\infty \left(\frac{\epsilon}{1-\delta}\right)^n \sum_{k=0}^{n-1} {n-1 \choose k} \theta^{k}\delta^{n-1-k}\\
                                                   &= (\delta + \theta)^2 \left(\frac{\epsilon}{1-\delta}\right)^2 \sum_{n=2}^\infty  \left(\frac{\epsilon}{1-\delta} (\delta+\theta)\right)^{n-2}\\
                                                   &= \left( \frac{\epsilon}{1-\delta}\right)^2 (\delta+\theta)^2 \frac{1-\delta}{1-\delta-\epsilon(\delta+\theta)},
  \end{align*}
  where the result holds only for $\epsilon(\delta+\theta)/\lvert 1-\delta\rvert < 1$.
\end{proof}

\begin{proof}[Proof of Theorem~\ref{t.tail perturbed ruin probability}]
  When \ptc\ has a phase-type distribution, then \eptc\ has also a phase-type distribution \cite{asmussen-RP}, and consequently it has an exponential decay rate. Thus, by the definition of the stationary excess claim sizes \emixc\ and Property~\ref{p.lighter tail subexponential}, we have
   \begin{equation}\label{e.tail perturbed stationary claim size distribution}
    \pr[(\emixc > u)] =\frac{(1-\epsilon)\delta}{(1-\epsilon)\delta+\epsilon\theta}\pr[(\eptc > u)] +\frac{\epsilon\theta}{(1-\epsilon)\delta+\epsilon\theta}\pr[(\ehtc > u)]
                \sim \frac{\epsilon\theta}{(1-\epsilon)\delta+\epsilon\theta}\pr[(\ehtc > u)],
  \end{equation}
  which implies by Property~\ref{p.subexponential class closed} that $\emixc \in \mathcal{S}$. When $\emixc \in \mathcal{S}$, it is known \cite{asmussen-RP} that
  \begin{equation}\label{e.tail perturbed ruin probability}
    \ruinmix \sim \frac{\rho_\epsilon}{1-\rho_\epsilon} \pr[(\emixc > u)],
  \end{equation}
  where $\rho_\epsilon=(1-\epsilon)\delta+\epsilon\theta<1$. Combining \eqref{e.tail perturbed stationary claim size distribution} and \eqref{e.tail perturbed ruin probability} yields the result.
\end{proof}

\begin{proof}[Proof of Theorem~\ref{t.tail discard}]
  The discard approximation \ruindis\ has a phase-type representation; therefore, it is of $o\left(\pr[(\ehtc > u)]\right)$. The same holds for the tail of the distribution of $\dissup[,0]+\dissup[,1]$. Moreover, since $\ehtc \in \mathcal{S}$, from Property~\ref{p.lighter tail subexponential} we obtain $\pr(\dissup[,0]+\dissup[,1]+\ehtc[1]> u) \sim \pr(\ehtc > u)$, which leads to the result by inserting these asymptotic estimates into \eqref{e.corrected discard}.
\end{proof}

\begin{proof}[Proof of Theorem~\ref{t.tail replace}]
  The class of phase-type distributions is closed under convolutions \cite{asmussen-RP}, which means that both $\repsup[,0]+\repsup[,1]$ and $\repsup[,0]+\repsup[,1]+\eptc[1]$ follow some phase-type distribution. Therefore, due to their exponential decay rate, \ruinrep, $\pr(\repsup[,0]+\repsup[,1]> u)$ and $\pr(\repsup[,0]+\repsup[,1]+\eptc[1]> u)$ are all of the order $o\left(\pr[(\ehtc > u)]\right)$. In addition, since $\ehtc \in \mathcal{S}$, we obtain from Property~\ref{p.lighter tail subexponential} that $\pr(\repsup[,0]+\repsup[,1]+\ehtc[1]> u) \sim \pr(\ehtc > u)$. Inserting these asymptotic estimates into \eqref{e.corected replace} leads to the result.
\end{proof}

\begin{proof}[Proof of Theorem~\ref{t.discard relative error}]

Let $X_{\epsilon,i}$ be an i.i.d.\ sequence such that $X_{\epsilon,i}\stackrel{d}{=} \dissup[,0] + \dissup[,i] +\ehtc[i]$, and similarly let $Y_{\epsilon,i}$ be an i.i.d.\ sequence such that $Y_{\epsilon,i}\stackrel{d}{=} \dissup[,0] +\ehtc[i]$.
Since \ehtc[i] is subexponential, and \dissup[,0] is light-tailed, according to Property~\ref{p.lighter tail subexponential}, $X_{\epsilon,i}, Y_{\epsilon,i}$ are subexponential as well. In order to prove Theorem~\ref{t.discard relative error}, we first need the following lemma.

\begin{lemma}\label{l.bounded subexponentials}
  There exists a constant $K_0$ independent of $\epsilon$, such that
  \begin{equation*}
    \frac{\pr(X_{\epsilon,1}+\dots +X_{\epsilon,n}>u)}{\pr(X_{\epsilon,1}>u)} \leq K_0^n,
  \end{equation*}
  for all $u$ and for all $n$.
\end{lemma}

\begin{proof}
  We follow a similar idea as the proof of Lemma 1.3.5 in \cite{embrechts-MEE}, which is not directly applicable, as $X_{\epsilon,i}$ depends on $\epsilon$. Let $F$ be the distribution function of $X_{\epsilon,i}$. We set $\alpha_n = \sup_u \disf{\com{\con{F}{}}}/\disf{\com{F}}$. Observe that
    \begin{align*}
        \frac{\disf{\com{\con[(n+1)]{F}{}}}}{\disf{\com{F}}} &= 1 + \int_0^u \frac{\disf[u-x]{\com{\con{F}{}}}}{\disf{\com{F}}} d\disf[x]{F}
                                                              = 1 + \int_0^u \frac{\disf[u-x]{\com{\con{F}{}}}}{\disf[u-x]{\com{F}}}\frac{\disf[u-x]{\com{F}}}{\disf{\com{F}}} d\disf[x]{F}\\
                                                             &\leq 1 + \frac{\alpha_n}{\disf{\com{F}}}\left( \disf{\com{\con[2]{F}{}}} - \disf{\com{F}} \right)
                                                              \leq 1+ \alpha_n (\alpha_2-1).
    \end{align*}
    Recursively, we find that
    \begin{equation*}
      \alpha_{n+1} \leq \sum_{k=0}^{n-2}(\alpha_2-1)^{k} + \alpha_2(\alpha_2-1)^{n-1}.
    \end{equation*}
    From Definition~\ref{d.subexponentiality}, we know that $\alpha_2 -1\geq 1$. So,
    \begin{align*}
      \alpha_{n+1} \leq \sum_{k=0}^{n-2} \alpha_2^{k} + \alpha_2^{n}
                   \leq \sum_{k=0}^{n} \alpha_2^{k} \leq \alpha_2^{n+1},
    \end{align*}
    therefore it suffices to show that $\alpha_2$ is bounded in $\epsilon>0$.

    To this end, observe that \dissup[,i] is stochastically decreasing in $\epsilon$ as it is the supremum of a compound Poisson process with arrival rate $\lambda(1-\epsilon)$. Therefore, the supremum that corresponds to the compound Poisson process with arrival rate $\lambda$ ($\epsilon=0$) is stochastically larger than all other suprema with $\epsilon>0$ and we denote it by \repsup[i]. Letting $S = \sum_{i=1}^4 \repsup[,i]$, we see that
    \begin{align*}
       \frac{\disf{\com{\con[2]{F}{}}}}{\disf{\com{F}}} &= \frac{\pr(X_{\epsilon,1}+X_{\epsilon,2}>u)}{\pr(X_{\epsilon,1}>u)}
                                                         \leq \frac{\pr(S + \ehtc[1]+\ehtc[2] > u)}{\pr(\ehtc[1] > u)}\\
                                                        &= \frac{\pr(S>u)}{\pr(\ehtc[1] > u)} + \int_0^u \frac{\pr(\ehtc[1]+\ehtc[2] > u-x)}{\pr(\ehtc[1] > u-x)}\frac{\pr(\ehtc[1] > u-x)}{\pr(\ehtc[1] > u)}d\pr(S\leq x) \\
                                                        &\leq \frac{\pr(S>u)}{\pr(\ehtc[1] > u)} + \underbrace{\sup_{u>0} \frac{\pr(\ehtc[1]+\ehtc[2] > u)}{\pr(\ehtc[1] > u)}}_{>1}
                                                                    \frac{1}{\pr(\ehtc[1] > u)} \int_0^u \pr(\ehtc[1] > u-x)d\pr(S\leq x)\\
                                                        &\leq \sup_{u>0} \frac{\pr(\ehtc[1]+\ehtc[2] > u)}{\pr(\ehtc[1] > u)} \sup_{u>0} \frac{\pr(S + \ehtc[1] > u)}{\pr(\ehtc[1] > u)}.
    \end{align*}
    Both suprema are finite since \ehtc[1] is subexponential and $S$ has a lighter tail than \ehtc[1]. This completes the proof of the lemma.
\end{proof}

We now proceed with the proof of Theorem~\ref{t.discard relative error}.
    Set $p_\epsilon = \frac{\epsilon\theta}{1-\delta+\epsilon\delta}$. Let $\eta$ be such that $p_\eta K_0=1/2$ and suppose $\epsilon<\eta$. Let $N$ be a random variable such that $\pr(N=n)=(1-p_\epsilon)p_\epsilon^n$. Observe that $\mixsup \stackrel{d}{=} \dissup[,0] + \sum_{i=1}^N Y_{\epsilon,i}$. For notational convenience, we assume that this equality holds almost surely through this proof.  This enables us to write
    \begin{equation*}
        \ruinmix - \cordis = \pr(\mixsup > u ; N\geq 2) - p_\epsilon^2 \pr(\dissup[,0]+Y_{\epsilon,1}>u),
    \end{equation*}
    so that
    \begin{align}
        \disf{R_{d,\epsilon}} &= \frac{\pr(\mixsup > u ; N\geq 2) - p_\epsilon^2 \pr(\dissup[,0] + Y_{\epsilon,1}>u)}{\pr(\mixsup > u)} \notag\\
                   &= \frac{\pr(\mixsup > u ; N\geq 2) - p_\epsilon^2 \pr(\dissup[,0] + Y_{\epsilon,1}>u)}{\pr(\mixsup > u ; N\geq 1)}\frac{\pr(\mixsup > u ; N\geq 1)}{\pr(\mixsup > u)}.\label{e.relative error}
    \end{align}
    Note that $\frac{\pr(\mixsup > u ; N\geq 1)}{\pr(\mixsup > u)}\leq 1$, where this ratio actually converges to $1$ as $u\rightarrow\infty$. To analyze the other fraction of \eqref{e.relative error}, the memoryless property of $N$ yields $\pr(\mixsup > u ; N\geq k) = p_\epsilon^k \pr(\mixsup +Y_{\epsilon,1}+\dots + Y_{\epsilon,k}> u)$ so
    \begin{align*}
        \frac{\pr(\mixsup > u ; N\geq 2)}{\pr(\mixsup > u ; N\geq 1)} &= p_\epsilon \frac{\pr(\mixsup + Y_{\epsilon,1}+Y_{\epsilon,2}> u)}{\pr(\mixsup + Y_{\epsilon,1} > u)}
                                                                             \leq p_\epsilon \frac{\pr(\mixsup + Y_{\epsilon,1} + Y_{\epsilon,2}+ > u)}{\pr(X_{\epsilon,1} > u)} \\
                                                                            &\leq p_\epsilon \pr(N=0) \frac{P(X_{\epsilon,1}+Y_{\epsilon,2} > u)}{\pr(X_{\epsilon,1} > u)}
                                                                            + p_\epsilon \sum_{n=1}^\infty \pr(N=n) \frac{\pr(X_{\epsilon,1}+\dots +X_{\epsilon,n+2}>u)}{P(X_{\epsilon,1}>u)}\\
                                                                            &\leq p_\epsilon \frac{\pr(X_{\epsilon,1}+Y_{\epsilon,2} > u)}{\pr(X_{\epsilon,1} > u)} + p_\epsilon^2 (1-p_\epsilon) K_0^3 \sum_{n=1}^\infty
                                                                            (p_\epsilon K_0)^{n-1} \\
                                                                            &\leq  p_\epsilon \frac{\pr(X_{\epsilon,1}+Y_{\epsilon,2} > u)}{\pr(X_{\epsilon,1} > u)} + p_\epsilon^2 2K_0^3.
    \end{align*}
    Finally, note that
    \begin{align*}
        \frac{p_\epsilon^2 \pr(\dissup[,0]+Y_{\epsilon,1}>u)}{\pr(\mixsup > u ; N\geq 1)} =  p_\epsilon \frac{\pr(M^\bullet_{\epsilon,0}+Y_{\epsilon,1}>u)}{\pr(M_\epsilon +Y_{\epsilon,1}> u )}.
    \end{align*}
    As before, we can show there exists a constant $K_1$ such that $\frac{\pr(\mixsup + Y_{\epsilon,1}> u )}{\pr(\dissup[,0]+Y_{\epsilon,1}>u)}\leq 1+ p_\epsilon K_1$. Putting everything together, we conclude that
    \begin{align*}
        \disf{R_{d,\epsilon}} &\leq p_\epsilon \frac{\pr(X_{\epsilon,1}+Y_{\epsilon,2} > u)}{\pr(X_{\epsilon,1} > u)} + p_\epsilon^2 2K_0^3 - p_\epsilon\frac{1}{1+ p_\epsilon K_1} \\
                   &\leq p_\epsilon \left(\frac{\pr(X_{\epsilon,1}+Y_{\epsilon,2} > u)}{\pr(X_{\epsilon,1} > u)} -1\right)+ p_\epsilon^2 K,
    \end{align*}
    for some constant $K$, completing the proof.
\end{proof}

\begin{proof}[Proof of Theorem~\ref{t.perturbed ruin probability example}]
  The Laplace transform of the ruin probability $\mathcal{L}\{\ruinmix\}$ satisfies the equation
  \begin{equation}\label{e.lt perturbed ruin probability example}
     \mathcal{L}\{\ruinmix\} = \frac{\rho_\epsilon}{s} \left(1 - \frac{(1-\rho_\epsilon) \ltegc}{1-\rho_\epsilon \ltegc}\right),
  \end{equation}
  where $\rho_\epsilon = \frac{\lambda}{\mu\nu}\big(\mu + \epsilon(\nu-\mu)\big)$, and
  \begin{align*}
    \ltegc &= \frac{1}{\e U_\epsilon} \big( (1-\epsilon) \e B \lteptc + \epsilon \e C \ltehtc\big)
            = \frac{\mu\nu}{\mu + \epsilon(\nu-\mu)} \left( (1-\epsilon)\frac{1}{\nu} \frac{\nu}{s+\nu} +  \epsilon \frac{1}{\mu} \frac{\mu}{(\mu+\sqrt{s})(1+\sqrt{s})} \right)\\
           &= \frac{\mu\nu}{\mu + \epsilon(\nu-\mu)} \cdot \frac{(1-\epsilon)(\mu+\sqrt{s})(1+\sqrt{s}) +  \epsilon (s+\nu)}{(s+\nu)(\mu+\sqrt{s})(1+\sqrt{s})}.
  \end{align*}
  If we set $\df[s]{\lt{w}}=(1-\rho_\epsilon) \ltegc/\left(1-\rho_\epsilon \ltegc\right)$, then with simple calculations we find that
  \begin{equation*}
    \df[s]{\lt{w}} = \frac{\mu\nu - \lambda \left(\mu + \epsilon(\nu-\mu)\right)}{\mu + \epsilon(\nu-\mu)} \cdot
                               \frac{(1-\epsilon)(\mu+\sqrt{s})(1+\sqrt{s}) +  \epsilon (s+\nu)}{(s+\nu)(\mu+\sqrt{s})(1+\sqrt{s}) - \lambda (1-\epsilon)(\mu+\sqrt{s})(1+\sqrt{s}) -\lambda \epsilon (s+\nu)}.
  \end{equation*}
  The denominator of $\df[s]{\lt{w}}$,
  \begin{equation*}
    \df[\sqrt{s}]{d} = s^2 +(\mu+1) s\sqrt{s} + (\mu+\nu-\lambda) s + (\mu+1)(\nu - \lambda+\lambda\epsilon)\sqrt{s} + \big(\mu(\nu-\lambda)+\lambda\epsilon(\mu-\nu)\big),
  \end{equation*}
  is a fourth degree polynomial with respect to $\sqrt{s}$. Let its roots be given by $-\nu_i(\epsilon)$, $i=1,\dots,4$, and let $\df[s]{n}$ denote the numerator of $\df[s]{\lt{w}}$. Then,
  \begin{equation}\label{e.rational fraction example}
    \frac{\df[s]{n}}{\df[s]{d}} = \sum_{i=1}^4 \frac{a_i}{\sqrt{s} + \nu_i(\epsilon)}.
  \end{equation}
  Finally, the coefficients $a_i$ are determined by the following equations
  \begin{equation*}
    a_i = \lim_{\sqrt{s} \rightarrow -\nu_i(\epsilon)}\frac{\df[s]{n}}{\df[s]{d}}\left(\sqrt{s} + \nu_i(\epsilon)\right), \qquad i=1,\dots,4.
  \end{equation*}
  For $s=0$, from \eqref{e.rational fraction example} we get
  \begin{equation*}
    0= \df[0]{n} - \df[0]{d} \sum_{i=1}^4 \frac{a_i}{\nu_i(\epsilon)}=
    \mu + \epsilon(\nu-\mu) - \Big(\mu\nu - \lambda \big(\mu + \epsilon(\nu-\mu)\big)\Big)\sum_{i=1}^4 \frac{a_i}{\nu_i(\epsilon)}.
  \end{equation*}
  Substituting everything in \eqref{e.lt perturbed ruin probability example}, we find
  \begin{align*}
      \mathcal{L}\{\ruinmix\} =& \frac{1}{s}\frac{\lambda}{\mu\nu}\big(\mu + \epsilon(\nu-\mu)\big)
                           \left(1 -  \frac{\mu\nu - \lambda \big(\mu + \epsilon(\nu-\mu)\big)}{\mu + \epsilon(\nu-\mu)}\sum_{i=1}^4 \frac{a_i}{\sqrt{s} + \nu_i(\epsilon)}\right)\\
                        =& \frac{\lambda}{\mu\nu} \left( \frac{\mu + \epsilon(\nu-\mu)}{s} - \Big(\mu\nu - \lambda \big(\mu + \epsilon(\nu-\mu)\big)\Big)\sum_{i=1}^4 \frac{a_i}{s\big(\sqrt{s} + \nu_i(\epsilon)\big)} \right)\\
                        =& \frac{\lambda}{\mu\nu} \left( \frac{\mu + \epsilon(\nu-\mu)}{s} - \Big(\mu\nu - \lambda \big(\mu + \epsilon(\nu-\mu)\big)\Big)\sum_{i=1}^4 \frac{a_i}{\nu_i(\epsilon)s}\right)\\
                         &+ \frac{\lambda}{\mu\nu} \Big(\mu\nu - \lambda \big(\mu + \epsilon(\nu-\mu)\big)\Big)\sum_{i=1}^4 \frac{a_i}{\nu_i(\epsilon)}\frac{1}{\big(\sqrt{s} + \nu_i(\epsilon)\big)\sqrt{s}}\\
                        =& \frac{\lambda}{\mu\nu} \Big(\mu\nu - \lambda \big(\mu + \epsilon(\nu-\mu)\big)\Big)\sum_{i=1}^4 \frac{a_i}{\nu_i(\epsilon)}\frac{1}{\big(\sqrt{s} + \nu_i(\epsilon)\big)\sqrt{s}},
  \end{align*}
  Laplace inversion to $\mathcal{L}\{\ruinmix\}$ gives,
    \begin{equation*}
    \df{\psi_\epsilon}= \frac{\lambda}{\mu\nu} \Big(\mu\nu - \lambda \big(\mu + \epsilon(\nu-\mu)\big)\Big)\sum_{i=1}^4 \frac{a_i}{\nu_i(\epsilon)} \df[\nu_i^2(\epsilon) u]{\zeta}.
  \end{equation*}
\end{proof}

\begin{proof}[Proof of Theorem~\ref{t.loss corrected discard error bounds}]
  Using that conditional probabilities are less than or equal to 1, an upper bound for the error of the approximation \cordisloss\ is found as
  \begin{align*}
    \pr \big(\df[t]{N_\epsilon^H} \geq 2 \big) \Bigg[ & \pr\Bigg( \left. \sum_{i=1}^{\df[t]{N_\epsilon^P}} \ptc[i] + \sum_{i=1}^{\df[t]{N_\epsilon^H}} \htc[i]  >x \right|  \df[t]{N_\epsilon^H} \geq 2 \Bigg)
    - \pr\Bigg(\sum_{i=1}^{\df[t]{N_\epsilon^P}} \ptc[i] + \htc >x \Bigg) \Bigg]\\
    &\leq \pr\Bigg( \left. \sum_{i=1}^{\df[t]{N_\epsilon^P}} \ptc[i] + \sum_{i=1}^{\df[t]{N_\epsilon^H}} \htc[i]  >x  \right|  \df[t]{N_\epsilon^H} \geq 2 \Bigg) \pr \big(N_\epsilon^H(t) \geq 2 \big)\\
    & \leq  \pr \big(\df[t]{N_\epsilon^H} \geq 2 \big)
    = \sum_{k=2}^\infty \frac{\big(\lambda t \big)^k}{k!}\epsilon^k e^{-\lambda \epsilon t}
    = \epsilon^2 (\lambda t)^2\sum_{k=2}^\infty \frac{\big(\lambda t \big)^{k-2}}{k!}\epsilon^{k-2} e^{-\lambda \epsilon t}\\
    &\leq \epsilon^2 (\lambda t)^2\sum_{k=2}^\infty \frac{\big(\lambda t \big)^{k-2}}{(k-2)!}\epsilon^{k-2} e^{-\lambda \epsilon t}
    = \epsilon^2 (\lambda t)^2.
  \end{align*}
  Using the obvious relation
  \begin{align*}
     \pr\Bigg( \left. \sum_{i=1}^{\df[t]{N_\epsilon^P}} \ptc[i] + \sum_{i=1}^{\df[t]{N_\epsilon^H}} \htc[i]  >x \right|  \df[t]{N_\epsilon^H} \geq 2 \Bigg)
     \geq
      \pr\Bigg( \sum_{i=1}^{\df[t]{N_\epsilon^P}} \ptc[i] + \htc[1] + \htc[2]  >x  \Bigg)
      \geq
       \pr\Bigg( \sum_{i=1}^{\df[t]{N_\epsilon^P}} \ptc[i] + \htc[1]  >x \Bigg),
  \end{align*}
  it is easy to verify that the error is non-negative, which completes the proof.
\end{proof}

\begin{proof}[Proof of Theorem~\ref{t.LST mixed poisson}]
  First, we define the Laplace-Stieltjes transform of $S_{\df[t]{N_\epsilon^P}} = \sum_{i=1}^{\df[t]{N_\epsilon^P}} \ptc[i]$ as
  \begin{align*}
    \hat{F}_{\df[t]{N_\epsilon^P}}(s) = \int_0^\infty e^{-s x} d\pr\big(S_{\df[t]{N_\epsilon^P}}\leq x\big)
                                      = \sum_{k=0}^\infty \pr\big(\df[t]{N_\epsilon^P} = k\big) \ltptc{k}
                                      = \text{exp}\{ -\lambda (1-\epsilon) t \big(1-\ltptc{}\big) \}.
  \end{align*}
  Consequently, the Laplace-Stieltjes transform of \cordisloss[T] satisfies
  \begin{align*}
  \mathcal{L}\{\cordisloss[T]\} =& \int_{x=0}^{\infty} e^{-s x} \int_{t=0}^\infty e^{-\lambda \epsilon t}  d\pr\Bigg(\sum_{i=1}^{\df[t]{N_\epsilon^P}} \ptc[i] >x  \Bigg) d\pr(T\leq t)\\
                      &+ \int_{x=0}^{\infty} e^{-s x} \int_{t=0}^\infty \big( 1-e^{-\lambda \epsilon t} \big) d\pr\Bigg(\sum_{i=1}^{\df[t]{N_\epsilon^P}} \ptc[i] + \htc >x \Bigg) d\pr(T\leq t)\\
                     =& \int_{t=0}^\infty e^{-\lambda \epsilon t} d\pr(T\leq t) \int_{x=0}^{\infty} e^{-s x} d\pr\Bigg(\sum_{i=1}^{\df[t]{N_\epsilon^P}} \ptc[i] >x  \Bigg) \\
                      &+ \int_{t=0}^\infty \big( 1-e^{-\lambda \epsilon t} \big) d\pr(T\leq t) \int_{x=0}^{\infty} e^{-s x} d\pr\Bigg(\sum_{i=1}^{\df[t]{N_\epsilon^P}} \ptc[i] + \htc >x \Bigg) \\
                     =& \int_{t=0}^\infty e^{-\lambda \epsilon t} \frac{1 - e^{ -\lambda (1-\epsilon) t \big(1-\ltptc{}\big)}}{s} d\pr(T\leq t)\\
                      &+ \int_{t=0}^\infty \big( 1-e^{-\lambda \epsilon t} \big)  \frac{1 - e^{ -\lambda (1-\epsilon) t \big(1-\ltptc{}\big)} \lthtc}{s} d\pr(T\leq t) \\
                     =& \frac1s -\frac{1 - \lthtc}s \int_{t=0}^\infty e^{ -\lambda \big(1- (1-\epsilon) \ltptc{}\big) t} d\pr(T\leq t)
                       - \frac{\lthtc}s \int_{t=0}^\infty e^{ -\lambda (1-\epsilon) \big(1-\ltptc{}\big) t} d\pr(T\leq t)\\
                     =&  \frac1s -\frac{1 - \lthtc}s \lttime[\lambda \big(1- (1-\epsilon) \ltptc{}\big)] - \frac{\lthtc}s \lttime[\lambda (1-\epsilon) \big(1-\ltptc{}\big)],
  \end{align*}
  which completes the proof.
\end{proof}

%\section*{References}
%
%\bibliographystyle{apt}
%\bibliography{D:/elena}

\end{document}